\documentclass[a4paper]{article}

\title{Long paths and cycles in subgraphs of the cube}
\author{Eoin Long\thanks{St. John's College, Cambridge, United Kingdom. E-mail:
E.P.Long@dpmms.cam.ac.uk. Research is supported by a Benefactor Scholarship from
St. John's College, Cambridge.}}
\date{}

 \usepackage{amsmath, amsfonts}
 \usepackage{amssymb}
 \usepackage{amsthm}
 \usepackage{mathrsfs}
 \usepackage{graphicx}
 \usepackage{color}
 \usepackage{dsfont}

\begin{document}

\newtheorem{thm}{Theorem}[section]
\newtheorem{dfn}[thm]{Definition}
\newtheorem{lem}[thm]{Lemma}
\newtheorem{cor}[thm]{Corollary}
\newtheorem{prop}[thm]{Proposition}
\newtheorem{conjecture}[thm]{Conjecture}
\newtheorem*{exam}{Examples} %[thm]

\theoremstyle{remark} \newtheorem*{remark}{Remark}

\newcommand{\N}{\mathbb{N}}
\newcommand{\Z}{\mathbb{Z}}
\newcommand{\Q}{\mathbb{Q}}
\newcommand{\R}{\mathbb{R}}
\newcommand{\E}{\mathbb{E}}
\newcommand{\T}{\mathbb{T}}
\newcommand{\C}{\mathbb{C}}
\newcommand{\B}{\mathcal{B}}
\newcommand{\w}{\tilde}
\newcommand{\D}{\mathcal{D}}
\newcommand{\X}{\mathbf{X}}
\newcommand{\Y}{\mathbf{Y}}
\newcommand{\A}{\mathscr{A}}

\maketitle

\begin{abstract}
{
Let $Q_n$ denote the graph of the $n$-dimensional cube with vertex set
$\{0,1\}^n$ in which two vertices are adjacent if they differ in exactly one
coordinate. Suppose $G$ is a subgraph of $Q_n$ with average degree at least $d$.
How long a path can we guarantee to find in $G$? 

Our aim in this paper is to show that $G$ must contain an exponentially long
path. In fact, we show that if $G$ has minimum degree at least $d$ then $G$ must
contain a path of length $2^d-1$. Note that this bound is tight, as shown by a
$d$-dimensional subcube of $Q_n$. We also obtain the slightly stronger result
that $G$ must contain a cycle of length at least $2^d$. 
}
\end{abstract}

\section{Introduction}

Given a graph $G$ of average degree at least $d$, a classical result of Dirac
\cite{dirac} guarantees a path of length $d$ in $G$. Moreover, this bound is
best possible as can be seen from $K_{d+1}$.

Inside the cube $Q_n$ can we improve this bound? That is, given a subgraph $G$
of $Q_n$ with average degree at least $d$, what is the length of the longest
path in $G$? The edge isoperimetric inequality for the cube (\cite{bernstein},
\cite{harper}, \cite{hart}, \cite{lindsey}, see \cite{com} for background) says
that any subgraph of average degree at least $d$ must have size at least $2^d$.
In light of this, the above linear bound seems very weak. A natural subgraph of
$Q_n$ with average degree at least $d$ is the $d$-dimensional cube $Q_d$, the
analogue of the complete graph in $Q_n$, which contains a path of length
$2^d-1$. Must the size of the longest path in $G$ also be exponential?

The main result of this paper answers this question in the affirmative.

\begin{thm}
{
\label{mainthm}
Every subgraph $G$ of $Q_n$ with minimum degree $d$ contains a path of length
$2^{d}-1$.
}
\end{thm}

Note that this is best possible as shown by a $d$-dimensional subcube of $Q_n$.
In fact, the proof of Theorem \ref{mainthm} shows that we can always find a
longer path in $G$ unless it is isomorphic to $Q_d$. Using the well known fact
that every graph with average degree at least $d$ contains a subgraph with
minimum degree at least $\frac{d}{2}$ we obtain the following corollary to
Theorem \ref{mainthm}.

\begin{cor}
{
\label{corollary} Every subgraph $G$ of $Q_n$ with average degree at least $d$
contains a path of length at least $2^{\frac{d}{2}}-1$.
}
\end{cor}

We do not know a tight bound for average degree $d$. We also obtain the
corresponding result for the length of the longest cycle in subgraphs of $Q_n$
with large minimum degree.

\begin{thm}
{
\label{cyclethm}
Every subgraph $G$ of $Q_n$ with minimum degree $d$ contains a cycle of length
at least $2^d$.
}
\end{thm}

In Section 2 we give an overview of the proofs of Theorems \ref{mainthm} and
\ref{cyclethm}. The theorems themselves are then proved in Sections 3-7.

In Section 8 we show that the lower bound from Theorems \ref{mainthm} and
\ref{cyclethm} also extends to subgraphs of the grid graph $\Z^n$ and the
discrete torus $C_k^n$, for all $k\geq 4$. We also give a generalization of
Theorems \ref{mainthm} and \ref{cyclethm} to general `product-type' graphs in 
the following form:

\begin{thm}
{
\label{extension} Let $k\in \N$. Suppose $G$ is a graph with minimum degree at least 
$d$ and that $G$ has the following property:
\begin{quotation}
{\noindent Given any two vertices $x,y\in G$, there is a partition of $V(G)$
into two sets $X$ and $Y$ 
with $x\in X$ and $y\in Y$ such that $d_{G[X]}(v)\geq d(v)-k$ for all $v\in X$
and $d_{G[Y]}(v)\geq d(v)-k$ for all $v\in Y$.}
\end{quotation}
Then $G$ contains a path of length at least $2^{\frac{d}{k+2}}$.
}
\end{thm}

In Section 8 we also give some consequences of this theorem and make some conjectures.

% Section 2: Overview

\section{Overview}

As in the statement of Theorem \ref{mainthm}, let $G$ be a subgraph of $Q_n$
with $\delta (G)\geq d$. We will view the vertices of $Q_n$ as elements of the
power set of $[n]$, $\mathcal{P}[n]$. 

A plausible approach to proving Theorem \ref{mainthm} is to split $G$ along some
direction $i$ to obtain two induced subgraphs $G_1$ and $G_2$ consisting of
those vertices of $G$ respectively containing and not containing $i$, for some
$i\in [n]$. Provided such a direction is chosen to ensure that $G_1, G_2\neq
\emptyset $, we have $\delta (G_1),\delta (G_2)\geq d-1$ and by induction on Theorem
\ref{mainthm} we have a path of length $2^{d-1}-1$ in each subgraph. If we could
join these two paths into one we would clearly be done. However, as Theorem
\ref{mainthm} provides no information on where these paths start or end, we can
not expect to be able to do this. 

This suggests that we strengthen Theorem \ref{mainthm} to guarantee an
exponentially long path between \emph{any} two vertices $x$ and $y$ of $G$. In
general this is not possible -- for example, consider the graph $G'$ obtained by
removing all but one edge $xy$ of direction $d+1$ from the ($d+1$)-dimensional
cube $Q_{d+1}$. 

However this graph is not $2$-connected. The following theorem says that this is
the only obstruction to such a strengthening.

% Section 2: Main Theorem

\begin{thm}
{
\label{conn}
Let $G$ be a $2$-connected subgraph of $Q_n$ and $a$ and $b$ be distinct
vertices of $G$. Suppose that $d_G(z) \geq d$ for all $z \in G-\lbrace
a,b\rbrace$. Then $G$ contains an $a-b$ path of length at least $2^d-2$.
Furthermore, unless $G$ is isomorphic to $Q_d$ with $a$ and $b$ at even Hamming
distance, $G$ contains an $a-b$ path of length at least $2^{d}-1$.
}
\end{thm}

Note that we do not assume that $a$ or $b$ have degree at least $d$ in Theorem
\ref{conn}. This slight weakening of the minimum degree condition will allow us
to use induction on various subgraphs of $G$ which would otherwise not be
available. 

Before continuing with the overview we make a small diversion to introduce some
definitions: these are standard (e.g. see \cite{mgt}).
\vspace{.2cm}

% Section 2: Definitions of blocks, cutvertices, etc. as in MGT

A subgraph $B$ of a graph $G$ is a \emph{block} of $G$ if $B$ is either a bridge
of $G$ or forms a maximal $2$-connected subgraph of $G$. By maximality,
$|B_1\cap B_2|\leq 1$ for any two blocks $B_1$ and $B_2$ of $G$ and $G-E(B)$
contains no $x-y$ path between distinct vertices $x,y$ in a block $B$. Therefore
if any two blocks intersect, their common vertex must be a cutvertex and
conversely every cutvertex lies in at least two blocks. Since every cycle is
$2$-connected and an edge is a bridge iff it does not lie in any cycle, every
graph $G$ decomposes uniquely into its blocks $B_1,\ldots, B_p$ in the sense
that:
\begin{equation*}
{
E(G)=\bigcup_{i=1}^p E(B_i) \mbox{ and } E(B_i)\cap E(B_j)=\emptyset \mbox{ if }
i\neq j.
}
\end{equation*}

Suppose now that $G$ is connected. Let $\mathcal{B}(G)$, the
\emph{block-cutvertex graph} of $G$, be the bipartite graph with bipartition
$(\mathcal{B},\mathcal{C})$ where $\mathcal{B}$ is the set of blocks of $G$,
$\mathcal{C}$ is the set of cutvertices of $G$ with $Bc$ an edge if $c\in B$.
For a connected graph $G$, $\mathcal{B}(G)$ is a tree. 

The leaves of this tree are all elements of $\mathcal{B}$ and are called
\emph{endblocks}. Given an endblock $E$ we will denote its unique cutvertex by
cutv$(E)$. Note that a graph $G$ has only one endblock iff it is $2$-connected.
\vspace{.2cm}

We now return to the overview of the proof of Theorem \ref{conn}.

% Splitting lemma showing how to find Ga and Gb

\begin{lem}
 {
\label{splitting}
Let $G$ be a connected subgraph of $Q_n$ with $a$ and $b$ distinct vertices of
$G$. Then there exists a partition of $G$ into two connected subgraphs $G_a$ and
$G_b$ such that $a\in G_a$, $b\in G_b$ and for all $v\in G_c$, $d_{G_c}(v)\geq
d_G(v)-1$, where $c\in \{a,b\}$.
} 
\end{lem}

\begin{proof}
 {
Picking $i\in [n]$ such that $a$ and $b$ differ in coordinate $i$
and forming $G_1$ and $G_2$ as before, we have $a\in G_1$ and $b\in G_2$. Let
$C_b$ be the connected component of $G_2$ containing $b$. Taking $G_a$ to be the
connected component of $G-C_b$ containing $a$ and $G_b=G-G_a$ we are done.
 }
\end{proof}

We will refer to $i$ in the above proof as the splitting direction 
for $G_a$ and $G_b$.

% Section 2: Overview of the proof of previous theorem.

A central observation in the proof of Theorem \ref{conn} is that, provided
$d\geq 3$, given any endblock $E$ of $G_a$ with $a\notin E$, by induction on
Theorem \ref{conn}, $E$ contains a path of length at least $2^{d-1}-2$ from
$\mbox{cutv}(E)$ to any $y\in E-\mbox{cutv}(E)$ -- $d\geq 3$ here guarantees 
$E$ is $2$-connected and not a bridge. Since $G$ is $2$-connected there must
exist $y\in E-\mbox{cutv}(E)$ with a neighbour in $G_b$. Thus these endblocks
guarantee `endblock paths' of length at least $2^{d-1}-1$ from a point in $G_a$
to one in $G_b$. If we could find a path from $a$ to $b$ containing at least two
such endblock paths we would almost have our path (it may still be short 
two or three vertices to give the $2^{d}-2$ or $2^d-1$ bound).

For ease of exposition we will prove the following weakening of Theorem
\ref{conn} first. It will allow the reader to focus on the main ideas in the
proof of Theorem \ref{conn} without some distracting details necessary to 
ensure that an $a-b$ path formed from endblock paths is not slightly too 
short.

% Section 2: Weakened version of main Theorem --- 2^{d-1} instead of 2^d - 1

\begin{thm}
{
\label{connmod}
Let $G$ be a $2$-connected subgraph of $Q_n$ and $a,b\in V(G)$. Suppose that
$d_G(z) \geq d$ for all $z \in V(G)-\lbrace a,b\rbrace$. Then $G$ contains an
$a-b$ path of length at least $2^{d-1}$. 
}
\end{thm}

% Section 2: Explanation of degree problem

Another technicality that arises in the proof of Theorem  \ref{conn}
and \ref{connmod} is the possibility that the only choice of a splitting direction 
$i$ for $G_a$ and $G_b$ in Lemma \ref{splitting} above, leaves $a$ with just one 
neighbour in $G_a$ or $b$ with just one neighbour in $G_b$. While all cases can 
be dealt with simultaneously, we felt for clarity's sake it was easier to first 
restrict attention to the case where a splitting direction $i$ exists for $G_a$ 
and $G_b$ in which $d_{G_a}(a)\geq 2$ and $d_{G_b}(b)\geq 2$. 
\\

% Section 2: Explanation of what is proved in different sections 

Theorem \ref{connmod} is proved in Sections 3-6. Sections 3-5 will focus on the
case where we can find a partition direction $i$, such that
$d_{G_a}(a)\geq 2$ and $d_{G_b}(b)\geq 2$. Section 3 will describe the
block-cutvertex decomposition structure of $G_a$ and $G_b$ in the absence of an
$a-b$ path of length $2^{d-1}$ formed by joining at least two endblock paths
together, and Section 4 describes how the endblocks of $G_a$ interact with those
of $G_b$. In Section 5 we show that if $G$ does not contain a path from $a$ to
$b$ containing at least two endblock paths then the conditions of Theorem
\ref{connmod} hold for a smaller subgraph of $G$. This allows for an inductive
step and completes the proof of Theorem \ref{connmod} in this case.

Section 6 will allow us, using a small modification of the argument from
Sections 3-5, to extend from the case $d_{G_a}(a)\geq 2$ and $d_{G_b}(b)\geq 2$
to the general case, proving Theorem \ref{connmod}. 

Finally in Section 7 we show how to adjust the approach in Sections 3-6 to
obtain the optimal bound of Theorem \ref{conn}. 

To close this section we show that Theorem \ref{conn} implies Theorem
\ref{cyclethm}.
\vspace{0.2cm}

\noindent \textit{Proof of Theorem \ref{cyclethm}:} Take an endblock $E$ in the
block-cutvertex decomposition of $G$. Clearly $E$ is $2$-connected and all
vertices in $E-\mbox{cutv}(E)$ have at least $d$ neighbours in $E$. Pick a
neighbour $v$ of $\mbox{cutv}(E)$ in $E$. Then by Theorem \ref{conn} $G$
contains a $\mbox{cutv}(E)-v$ path $P$ of length at least $2^d-1$. Combining $P$
with the edge $\mbox{cutv}(E)v$ we obtain the desired cycle. \hspace{5cm}
$\square$

\section{Endblocks in $G_a$ and $G_b$}

To begin we introduce some useful definitions.

\begin{dfn}
{
\emph{
Let $E$ be an endblock in the block-cutvertex decomposition of $G_a$ ($G_b$).
The \emph{interior} of $E$ is the set $\mbox{int}(E)=E-\mbox{cutv}(E)$. A vertex
$x\in \mbox{int}(E)$ is said to be an \emph{exit vertex} of $E$ if $x$ has a
neighbour in $G_b$ ($G_a$). If this neighbour exists, it is unique and is
denoted by $p(x)$, $x$'s \emph{partner}.}
}
\end{dfn}

\begin{dfn}
{
\label{why}
\emph
{
Body$(a)$ is the intersection of all blocks of $G_a$ containing $a$. Let $\mbox
{Core}(a)$ consist of those vertices in $\mbox {Body}(a)$ that are not
cutvertices of $G_a$.
}
}
\end{dfn}

\begin{dfn}
{
\emph
{
A subgraph $K$ of $G_a$ is said to be a \emph{limb} of $a$ if:}
\begin{itemize}
\item \emph{$a$ is a cutvertex of $G_a$ and $K=G\lbrack C\cup \lbrace a\rbrace
\rbrack $ where $C$ is a connected component of $G_a-a$}

\item \emph{$a$ is not a cutvertex of $G_a$ and $K=G\lbrack C \rbrack $ where
$C$ is a connected component of $G_a-\mbox {Core}(a)$.}
\end{itemize}
\emph{
The \emph{joint} of a limb $K$, Joint$(K)$, is the unique vertex $v\in K\cap
\mbox {Body}(a)$.
}
}
\end{dfn}

\begin{figure}
\centering
\input{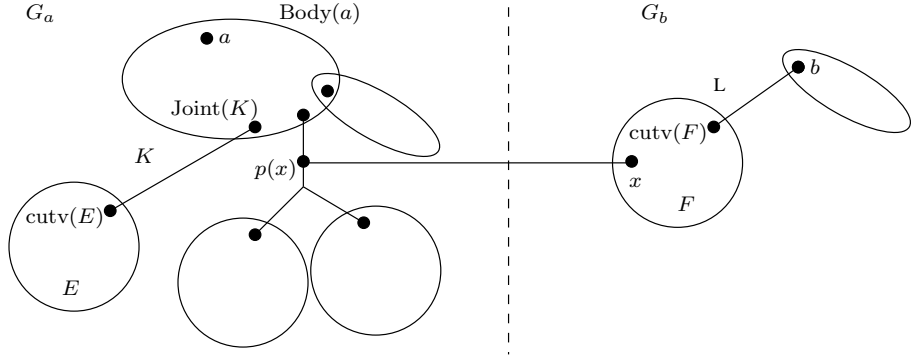}
\caption{The diagram displays various parts of $G_a$ and $G_b$. The broken line
separates $G_a$ and $G_b$. In $G_a$, $\mbox{Body}(a)\neq \{a\}$ and $a$ has
three limbs. In $G_b$, $b$ is a cutvertex and one of its limbs $L$ contains an
endblock $F$ with exit vertex $x$.\label{intro}}
\end{figure}

The reader may find it helpful to examine Figure \ref{intro}. The circles and
ellipses will always denote blocks in the block-cutvertex decomposition of the graph.
\vspace{0.2cm}

The proof of Theorem \ref{connmod} will proceed by induction on $d$. The case $d=2$ 
follows from Menger's theorem, as if $G$ is $2$-connected it
contains two disjoint $a-b$ paths, one of which must have length at least $2$.
We will suppose for contradiction that the Theorem fails for some $d>2$ and take 
$G$ to be a minimal counterexample so
that Theorem \ref{connmod} holds for all smaller degrees and all graphs $G'$
with $|G'|<|G|$. The following lemma will be the main step in the proof of Theorem 
\ref{connmod}. Its proof will be the aim of the next three sections.

% New inductive lemma

\begin{lem}
\label{splittingdirection}
Let $G$ be a $2$-connected subgraph of $Q_n$ and $a,b\in V(G)$ such that $d(v)\geq d$ for all 
$v\in V(G)-\{a,b \}$, where $d\geq 3$. Suppose that Theorem \ref{connmod} is true for smaller 
degrees and all graphs $G'$ with $|G'|<|G|$. Suppose furthermore that there exists a splitting 
direction $i$ for $G_a$ and $G_b$ in Lemma \ref{splitting} for which $d_{G_a}(a)\geq 2$ and 
$d_{G_b}(b)\geq 2$. Then $G$ contains an $a-b$ path of length at least $2^{d-1}$.
\end{lem}

Note that it follows from Lemma \ref{splittingdirection} that if $d_G(a)\geq 3$ and $d_G(b)\geq 3$, 
$G$ contains an $a-b$ path of length at least $2^{d-1}$. Indeed, taking any direction $i$ on 
which $a$ and $b$ differ as the splitting direction in the proof of Lemma \ref{splitting}, we have 
$d_{G_a}(a)\geq 2$ and $d_{G_b}(b)\geq 2$. Lemma \ref{splittingdirection} therefore applies and 
gives an $a-b$ path of length at least $2^{d-1}$, as claimed.

Over the next three sections we will establish some results which will be used in the proof of Lemma 
\ref{splittingdirection} in Section 5. The first of these describes the block structure 
of $G_a$ provided we cannot use endblock paths to form an $a-b$ path of length at least 
$2^{d-1}$.

\begin{lem}
{
\label{list}
Suppose that $G, a, b, G_a$ and $G_b$ are as in the statement of Lemma \ref{splittingdirection}. 
If $G$ does not contain an $a-b$ path of length at least $2^{d-1}$ then the following hold:
\begin{itemize}
{
\item[\emph{(i)}] Every endblock of $G_a$ which does not contain $a$ in its
interior must contain at least two exit vertices.

\item[\emph{(ii)}] $G_a$ is not $2$-connected.

\item[\emph{(iii)}] $a$ does not lie in the interior of an endblock in $G_a$.

\item[\emph{(iv)}] $a$ must have at least two limbs.
}
\end{itemize}
}
\end{lem}

\begin{proof}
{

(i)
Suppose not and let $E$ be such an endblock. By the $2$-connectivity of $G$, $E$
must contain an exit vertex $x$. If $x$ is its only exit vertex then every
$v\in E-\lbrace {\mbox{cutv}(E),x}\rbrace $ has degree least $d$ in $E$ --
such $v$ must exist since $d\geq 3$. Since $|E|<|G|$, $E$ contains a path 
$P_2$ of length at least $2^{d-1}$ from $\mbox{cutv}(E)$ to $x$. Joining $a$ to 
$\mbox{cutv}(E)$ in $G_a$ by a path $P_1$ and $p(x)$ to $b$ in $G_b$ by a path 
$P_3$ we have created a path $P_1P_2P_3$ of length at least $2^{d-1}$ from $a$ 
to $b$, a contradiction.

(ii)
Suppose $G_a$ is $2$-connected. First consider the case where $G_b$ is not
$2$-connected. Let $E$ be an endblock in $G_b$ not containing $b$ in its
interior and take $x$ to be an exit vertex of $E$ with $p(x)\neq a$ -- this
exists by (i). Since Theorem \ref{connmod} holds for $d-1$, there 
are paths $P_1$ in $G_a$ from $a$ to $p(x)$ and $P_2$ in $E$ from $x$ 
to $\mbox{cutv}(E)$ both of length at least $2^{d-2}$. Taking a path $P_3$ from 
$\mbox{cutv}(E)$ to $b$ in $G_b$ we have constructed a path $P=P_1p(x)xP_2P_3$ 
from $a$ to $b$ of length at least $2^{d-1}$, a contradiction. 

If $G_b$ is $2$-connected, then the same proof as in (i) shows that $G_b$ must
contain two exit vertices, one of which, $x$, has $x\neq b$ and $p(x) \neq a$.
Again as Theorem \ref{connmod} holds for $d-1$, we obtain endblock paths 
from $a$ to $p(x)$ in $G_a$ and from $x$ to $b$ in $G_b$ both of length at least 
$2^{d-2}$. Joining the two with edge $xp(x)$, $G$ again contains an $a-b$ path of 
length at least $2^{d-1}$, a contradiction.

(iii)
Suppose $a$ lies in the interior of an endblock $E$ of $G_a$. As
$d_{G_a}(a)\geq 2$ $E$ is $2$-connected. As Theorem \ref{connmod} holds for 
$d-1$, we have an endblock path $P_1$ from $a$ to $\mbox{cutv}(E)$ in $E$ 
of length at least $2^{d-2}$. From (ii) $G_a$ is not $2$-connected and so 
it contains a second endblock $E'$, with an exit vertex $x$. Again since Theorem 
\ref{connmod} holds for $d-1$, $E'$ contains an endblock path $P_3$ from 
$\mbox{cutv}(E')$ to $x$ of length $2^{d-2}$. Join $\mbox{cutv}(E)$ to 
$\mbox{cutv}(E')$ by a path $P_2$ in $G_a$ and $p(x)$ to $b$ by a path $P_4$ in 
$G_b$. Combining all of these paths we have a path $P_1P_2P_3xp(x)P_4$ from $a$ 
to $b$ of length at least $2^{d-1}$, a contradiction.

(iv)
This follows from (ii) and (iii) as if $G_a$ is not $2$-connected and $a$ does
not lie in the interior of any endblock, $a$ must have at least two limbs.
}
\end{proof}

\begin{figure}
\centering
\input{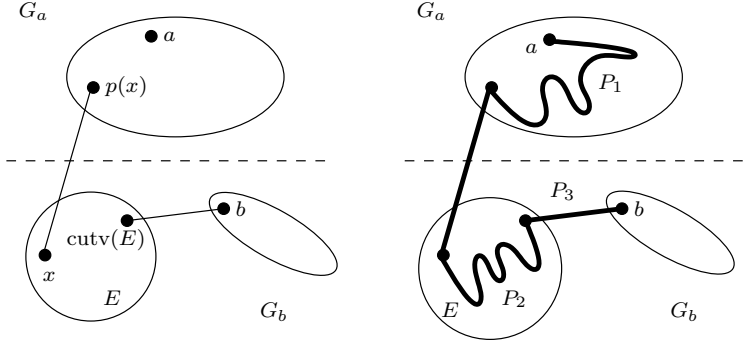}
\caption{Path $P$ constructed in Lemma \ref{list}(ii). Curved paths like $P_1$
and $P_2$ will represent endblock paths of length at least $2^{d-2}$
throughout.}
\end{figure}

Note that by symmetry of $a$ and $b$, Lemma \ref{list} also applies on replacing
$a$ with $b$. The next proposition gives a simple case  in which we can use
endblock paths to build our path of length $2^{d-1}$ from $a$ to $b$.

\begin{prop} 
{
\label{int}
Let $G, a, b, G_a$ and $G_b$ be as in the statement of Lemma \ref{splittingdirection}. Suppose 
$G$ does not contain an $a-b$ path of length at least $2^{d-1}$. Then for any exit vertex $x$ 
of an endblock $E$ of $G_a$, $p(x)$ can never lie in the interior of an endblock $F$ of $G_b$.
}
\end{prop}

\begin{proof}
{
From Lemma \ref{list}(iii) $a\notin \mbox{int}(E)$ and $b\notin \mbox{int}(F)$.
Pick a path $P_1$ in $G_a$ from $a$ to $\mbox{cutv}(E)$ and a path $P_4$ in
$G_b$ from $\mbox{cutv}(F)$ to $b$. Since $E$ is $2$-connected and all $v\in
E-\lbrace \mbox{cutv}(E),x\rbrace $ have degree at least $d-1$ in $G[E]$, 
Theorem \ref{connmod} gives a path $P_2$ of length at least $2^{d-2}$ from
$\mbox{cutv}(E)$ to $x$. Similarly $F$ contains a path $P_3$ of length at
least $2^{d-2}$ from $p(x)$ to $\mbox{cutv}(F)$. Combining these gives an $a-b$
path $P=P_1P_2xp(x)P_3P_4$ of length at least $2^{d-1}$, a contradiction.
}
\end{proof}

\section{The Interaction Digraph}

Throughout this section, $G$ will be a $2$-connected subgraph of $Q_n$ containing 
vertices $a$ and $b$, with $d_{G_a}(a)\geq 2$, $d_{G_a}(a)\geq 2$ and $d(v)\geq d$ 
for all $v\in V(G)-\{ a,b\} $. We will also assume that Theorem \ref{connmod} holds for 
all smaller degrees and for all graphs $G'$ with $|G'|<|G|$.

Let $K_1,\ldots ,K_r$ be the limbs of $a$ and $L_1,\ldots ,L_s$ be the limbs of
$b$. Lemma \ref{list}(iv) shows that $r,s\geq 2$.

We form an auxiliary bipartite multidigraph $H=(A,B,\overrightarrow {E})$ which
will represent the interaction between the limbs and cores of $a$ and $b$. Let
$A=\lbrace {K_1,\ldots ,K_r }\rbrace $ and $B=\lbrace {L_1,\ldots ,L_s }\rbrace
$. Additionally, adjoin $\mbox {Core}(a)$ to $A$ and $\mbox {Core}(b)$ to $B$ if
they are non-empty. Given an endblock $E$ of $G_a$, there exists an exit vertex
$x$ with $x\neq a$ and $p(x)\neq b$ by Lemma \ref{list}(i) and (iii). Pick
exactly one such exit vertex $x_E$ for each such endblock $E$ and adjoin a
directed edge to $H$ from $K$ to $W\in B$ where $E$ is contained in limb $K$ and
$p(x_E)\in W$. Similarly, for each endblock $F$ in $L$ we pick an exit vertex
$y_F\in F$ with $p(y_F) \neq a$ and add a directed edge to $H$ from $L$ to $V$
where $p(y_F)\in V$. 

Note that by Proposition \ref{int} we never choose an exit vertex $x_E$ for some
$E$ and $y_F$ for some $F$ such that $p(x_E)=y_F$. Also, since any limb of $a$ or
$b$ contains an endblock, every limb vertex in $H$ must have outdegree at least
one and core vertices have no outneighbours. 

We shall study the component structure of $H$. The next two lemmas say that this
is very restricted. Together they will allow us to find a connected
component $C$ of $H$ consisting entirely of limbs. The inductive step in Section
5 will take place on the subgraph of $G$ corresponding to this $C$.

As $H$ is a multidigraph, we stress that in the next lemma, by a path we
mean a path without repeated vertices.

\begin{lem}
{
\label{paths3}
Let $G, a, b, G_a$ and $G_b$ be as above. Suppose $G$ does not contain an $a-b$ path of length  
at least $2^{d-1}$. Then $H$ does not contain an undirected path of length three.
}
\end{lem}

\begin{proof}
{
Suppose for contradiction that we have such a path $Q=V_0V_1V_2V_3$ in $H$ and assume
$V_0\in A$. Each directed edge $\overrightarrow {VW}$ of $Q$ gives an endblock
in $V$ with exit vertex $x$, such that $p(x)\neq b$ and $p(x)\in W$. These
endblocks are distinct by the construction of $H$ and since Theorem \ref{connmod} 
holds for $d-1$, in each we can find an endblock path of length at least $2^{d-2}$ 
from its cutvertex to this exit vertex. We claim that we can form an $a-b$ path $P$ 
which extends all three of these paths. As such a path has length at least 
$3(2^{d-2})>2^{d-1}$, this contradicts the hypothesis and proves the lemma.

We will construct our path by forming paths $P_i$ in each $V_i$ and eventually
join them into one. The start point of $P_i$ will be denoted by $a_i$ and its
end point by $b_i$. We first choose these vertices.

If $\overrightarrow {V_iV_{i+1}}$ is an edge of $Q$ there is an endblock $E$ in
$V_i$ with an exit vertex $x_E$ such that $p(x_E)\in V_{i+1}$. In this case let
$b_i=x_E$ and $a_{i+1}=p(x_E)$. If $\overleftarrow {V_iV_{i+1}}$ is an edge of $Q$
this gives an endblock $E$ in $V_{i+1}$ with an exit vertex $x_E$ such that
$p(x_E)\in V_{i}$. In this case let $b_i=p(x_E)$ and $a_{i+1}=x_E$. We set
\[a_0 = \left\{ \begin{array}{ll}
         \mbox{Joint$(V_0)$} & \mbox{if $V_0$ is a limb of $a$};\\
        b_0 & \mbox{if $V_0=$Core$(a)$}\end{array} \right. \]
\[b_3 = \left\{ \begin{array}{ll}
         \mbox{Joint$(V_3)$} & \mbox{if $V_3$ is a limb of $b$};\\
        a_3 & \mbox{if $V_3=$Core$(b)$}.\end{array} \right. \] 

\noindent Note that $b_i$ and $a_{i+1}$ are adjacent for $i\in \lbrace
0,1,2\rbrace $ and $a,b\notin \{b_0,a_1,b_1,a_2,b_2,a_3\}$. 

We now build the paths $P_i$ from $a_i$ to $b_i$ in each $V_i$, where $V_i$ is a
limb. We claim we can choose $P_i$ so that neither $a$ nor $b$ are interior
vertices of $P_i$ (that is, they can lie on $P_i$, but only as end vertices) such 
that
$P_i$ has length at least $2^{d-2}$ if $V_i$ has one outneighbour on $Q$ and
$2^{d-1}$ if $V_i$ has two outneighbours on $Q$. Indeed, if $V_i$ has exactly one 
outneighbour in $Q$ then
exactly one of $a_i$ or $b_i$ must be an exit vertex of an endblock $E$ of $V_i$.
Without loss of generality this is $a_i$. We must also have $b_i\notin \mbox{int}(E)$. 
Indeed, by definition $b_3$ never lies in the interior of an endblock, so $i\leq 2$ 
and $a_{i+1}$ must be an exit vertex for a endblock in $V_{i+1}$. But as $b_i$ and 
$a_{i+1}$ are adjacent, this contradicts Propostion \ref{int}. Therefore, since 
Theorem \ref{connmod} holds for $d-1$, $E$ contains a path of length
$2^{d-2}$ from $a_i$ to the cutv$(E)$. Since $V_i-\lbrace a,b\rbrace $ is
connected for all $i$ from the definition of a limb, we can extend this path
from cutv$(E)$ to $b_i$ in $V_i$ as required. The case where $V_i$ has two 
outneighbours in $Q$ is identical, using the same argument in two endblocks of 
$V_i$ and joining their cutvertices in $V_i$.

Finally we combine the $P_i$ paths. We first deal with the case where neither
Core$(a)$ nor Core$(b)$ occur as interior vertices of $Q$. Combining the paths
above we have an $a_0-b_3$ path $P'=P_0b_0a_1P_1b_1a_2P_2b_2a_3P_3$. If
$\mbox{Body}(a)=\lbrace a\rbrace $ then $P'$ starts at $a$ so we only need to
extend $P'$ to start at $a$ when 
$\mbox{Body}(a)\neq \lbrace a\rbrace $. In $P'$ as constructed above, $\mbox
{Body}(a)\cap P'$ contains $a_0$ and at most one other vertex - indeed as the paths $P_i$ 
above always lie entirely inside $V_i$, they can only intersect $\mbox{Body}(a)$ 
in $\mbox{Joint}(V_i)$ and therefore $P'$ contains at most $a_0$ and $\mbox{Joint}(V_2)$. 
Since $\mbox{Body}(a)$ is $2$-connected it contains a path $P_1'$ from $a$ to $a_0$
avoiding $\mbox{Joint}(V_2)$. Finding a similar path $P_2'$ from $b_3$ to $b$ in 
$\mbox{Body}(b)$ if $\mbox {Body}(b)\neq \lbrace b\rbrace $ we may take
$P=P_1'P'P_2'$.

\begin{figure}
	\centering
	\input{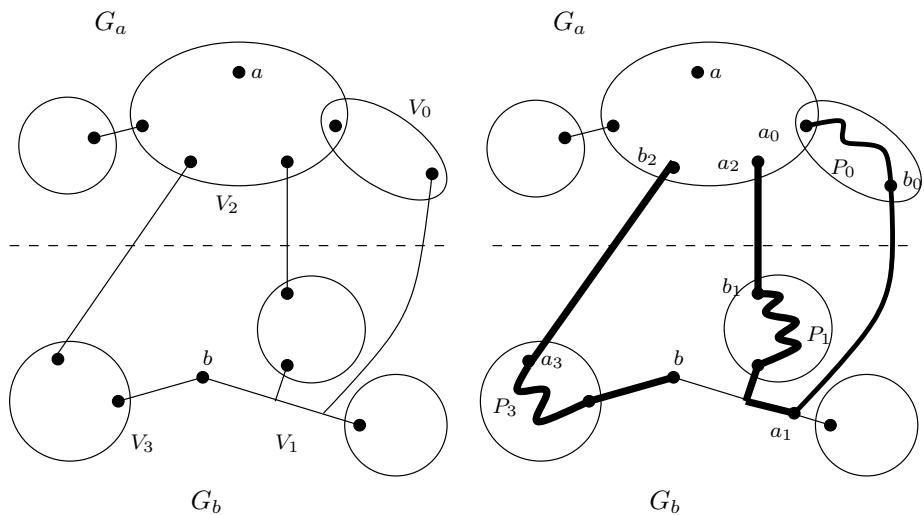}
	\caption{An illustration of Lemma \ref{paths3} in the case where
$V_2=\mbox{Core}(a)$ and $V_0V_1$, $V_1V_2$ and $V_3V_2$ are directed edges of
$Q$. As in the proof of Lemma \ref{paths3}, 2-connectivity can be used in
$\mbox{Body}(a)$ to find vertex disjoint paths from $\{a_0,a_2\}$ to
$\{b,b_2\}$.}
\label{picpaths3}
\end{figure}

If $Q$ contains one of the Core vertices, without loss of generality let it be
Core$(a)$. If Core$(a)$ occurs as an interior vertex of $Q$, it must be $V_2$.
Body$(a)$ then contains distinct $a_0,a_2,b_2$ and we have two paths
$P_1'=P_0b_0a_1P_1b_1a_2$ from $a_0$ to $a_2$ and $P_2'=b_2a_3P_3$ from $b_2$ to
$b_3$ as in Figure \ref{picpaths3}. From the choice of the $a_2$ and $b_2$ above
and the fact that $a$ is not a cutvertex we have $a\notin \lbrace
a_0,a_2,b_2\rbrace $. Therefore by 2-connectivity Body$(a)$ contains two vertex
disjoint paths from $\lbrace a_0,a_2\rbrace $ to $\lbrace a,b_2\rbrace $.
Piecing these paths together with $P_1'$ and $P_2'$ we obtain an $ab_3$-path
$P'$. If Body$(b)=\lbrace b\rbrace $ we are done since $b=b_3$. Otherwise we
extend $P'$ using 2-connectivity as above to find an $a-b$ path of length at
least $2^{d-1}$, contradicting the choice of $G$.
}
\end{proof}

Note that Lemma \ref{paths3} guarantees that $H$ has at least two connected
components. The next lemma further limits $H$. Its proof is very similar to that
of Lemma \ref{paths3}.

\begin{lem} 
{
\label{body}
Let $G, a, b, G_a$ and $G_b$ be as above. Suppose that $G$ does not contain an $a-b$ path 
of length at least $2^{d-1}$. Furthermore, suppose that $\mbox{Body}(a)\neq \lbrace a\rbrace $. 
Then no component of $H$ contains two vertices of $A$.
}
\end{lem}

\begin{proof}
{
Suppose $H$ has such a component $C$. Then, since $H$ does not contain a path of length 
three by Lemma \ref{paths3}, $C$ consists
of vertices $V_1,\ldots ,V_t$ in $A$ and a vertex $W$ in $B$. At most one of
$V_1,\ldots ,V_t, W$ can be a core vertex as there is no edge between Core$(a)$
and Core$(b)$ in $H$. 

If $W=\mbox{Core}(b)$ then $V_1$ and $V_2$ must be limbs and these guarantee two
vertex disjoint paths $P_1$, $P_2$ from vertices $a_1,a_2\in \mbox{Body}(a)$ to
vertices $b_1,b_2\in \mbox{Body}(b)$ both of length at least $2^{d-2}$ with 
$|P_i\cap \mbox {Body}(c)|=1$ for $i=1,2$ and $c\in \{a,b\}$. As $b$ has 
at least two limbs by Lemma \ref{list}(iv) and by Lemma
\ref{paths3} $C$ cannot contain both of these, $H$ must contain a second component $C'$
containing a limb of $b$. This guarantees the existence of a third path $P_3$
from a vertex $a_3\in \mbox{Body}(a)$ to $b_3\in \mbox{Body}(b)$ of length
$2^{d-2}$ again with $|P_3\cap \mbox {Body}(c)|=1$ for $c\in \{a,b\}$ which is disjoint 
from $P_1$ and $P_2$. Using identical 2-connectivity arguments 
in both $\mbox{Body}(a)$ and $\mbox{Body}(b)$ as in Lemma \ref{paths3} we can combine these 
three paths into one from $a$ to $b$, contradicting the hypothesis.

If $W\neq \mbox{Core}(b)$ then $C$ guarantees a path $P_1$ of length $2^{d-1}$
between two vertices $a_1$ and $a_2$ in $\mbox{Body}(a)$ with $|P_1\cap \mbox{Body}(a)|=2$, 
$b\notin P_1\cap \mbox{Body}(b)$ and  $|P_1\cap \mbox{Body}(b)|\leq 1$. Again from a second
connected component of $H$ we obtain a disjoint path $P_2$ from an element
$a_3\in \mbox{Body}(a)$ to $b_1\in \mbox{Body}(b)$ with $|P_2\cap \mbox{Body}(c)|\leq 1$ 
for $c\in \{a,b\}$. Once more, with an
application of 2-connectivity in $\mbox{Body}(a)$ and a possible application in
$\mbox{Body}(b)$ we find an $a-b$ path extending both $P_1$ and $P_2$, a
contradiction.
}
\end{proof}

Again the same applies switching $a$ with $b$. As mentioned before Lemma
\ref{paths3} the previous two lemmas imply that $H$ contains a connected
component $C$ consisting entirely of limbs.

\begin{cor} 
{
\label{limbcomponent}
Let $G, a, b, G_a$ and $G_b$ be as above. Suppose that $G$ does not 
contain an $a-b$ path of length $2^{d-1}$. Then the interaction digraph $H$ of 
$G$ has at least two connected components, one of which $C$ consists entirely 
of limbs.
}
\end{cor}

\begin{proof}
Since $|A|,|B|\geq 2$, if $H$ is connected it contains an undirected path of 
length three, contradicting Lemma \ref{paths3}. Therefore $H$ has at least 
two connected components, as claimed.
If $H$ does not contain a component consisting entirely of limbs, each component 
of $H$ contains one of Core$(a)$ or Core$(b)$. But then $H$ has exactly 
two connected components, one containing Core$(a)$ and one containing Core$(b)$. 
But as $A$ contains Core$(a)$ and at least two limbs, two of these must lie in the 
same connected component of $H$, contradicting Lemma \ref{body}.
\end{proof}

We will write $G_C$ for the subgraph $G[\cup _{W\in C}V(W)]$ of $G$. We note that $G_{C}$ 
must contain exactly one vertex $a_C$ in $\mbox{Body}(a)$ and one vertex $b_C$ in 
$\mbox{Body}(b)$ -- if Body$(a)= \lbrace a\rbrace$ then $a_C=a$, if not then by Lemma 
\ref{body} $A\cap C=\lbrace V\rbrace $ and we may take $a_C=\mbox{Joint}(V)$.

% Section 5: The Inductive Step

\section {The Inductive Step}

% Motivation for Inductive Lemma

Suppose that $G, a, b, G_a$ and $G_b$ satisfy the hypothesis of Lemma \ref{splittingdirection} 
but $G$ does not contain a path of length $2^{d-1}$. Then we may apply Corollary 
\ref{limbcomponent} to find a component $C$ of $H$ consisting entirely of limbs. 
Our final lemma before the proof of Lemma \ref{splittingdirection} finds a subgraph 
of $G_C$ which either also satisfies the conditions of Theorem \ref{connmod} or builds 
half of the $a-b$ path we are looking for from any edge entering it. Before 
stating it we give one last definition.

% Definition 5.1: span _G(S)

\begin{dfn}
{
\emph{Given a graph $G$ and $S\subset V(G)$ define the span$_{G}(S)$ to be
the subset of $V(G)$ consisting of all vertices which lie on a path between two
elements of $S$.}
}
\end{dfn}

Note that we include paths of length zero in this definition, so that 
$S\subset \mbox{span}_G(S)$.

% Lemma 5.2: The Inductive Lemma

\begin{lem}
{
\label{nice}
Let $G$ be a $2$-connected subgraph of $Q_n$ containing vertices $a$ and $b$ such 
that $d_{G_a}(a)\geq 2$, $d_{G_b}(b)\geq 2$ and $d(v)\geq d$ for all 
$v\in V(G)-\{ a,b\} $. Suppose that Theorem \ref{connmod} holds for all smaller 
degrees and for all graphs $G'$ with $|G'|<|G|$. Suppose furthermore that $G$ 
does not contain an $a-b$ path of length at least $2^{d-1}$. Then taking $C$ 
as in Corollary \ref{limbcomponent}, $G_C$ has a $2$-connected subgraph $J$ 
containing two vertices $a'\in G_a$ and $b'\in G_b$ with the following properties:

% Two properties of subgraph J

\begin{itemize}
{
\item[\emph{(i)}] every vertex $v\in J - \lbrace a',b'\rbrace$ has degree at
least $d-1$ in $J$ and all the neighbours of $v$ in $G_C$ are contained in $J$.

\item[\emph{(ii)}] for any vertex $v\in J-\lbrace a',b'\rbrace $, $J$ contains
an $a'-v$ path not containing $b$ and a $b'-v$ path not containing $a$, both of
length at least $2^{d-2}$.
}
\end{itemize}
}
\end{lem}

% Proof of the Inductive Lemma (5.2)

\begin{proof}
{
From Lemma \ref{paths3}, $C$ cannot contain two limbs of both $a$ and $b$. 
We may therefore assume that $C$ consists of one limb $K$ of $a$ and limbs 
$L_1,\ldots ,L_t$ of $b$. 

For each $i\in [t]$ we define $S_i\subset K$ and $T_i\subset L_i$ as follows:

\begin{equation*}
 S_i := \{v\in K: v \mbox{ has a neighbour }p(v)\in L_i-\{b\}\}
\end{equation*}
and 
\begin{equation*}
 T_i := \{w\in L_i: w \mbox{ has a neighbour }p(w)\in K-\{a\}\}.
\end{equation*}

Now each limb in the interaction digraph has at least one outneighbour. We claim that 
for each endblock $E\in K$ there exists some $i\in [t]$ with $|S_i|\geq 2$ such that 
$\mbox{int}(E)\cap S_i \neq \emptyset $. Indeed, from construction of the interaction 
graph, $E$ contributes a directed edge from $K$ to $L_i$ for some $i\in [t]$. This gives 
an exit vertex $x_E\in \mbox{int}(E)$ with $p(x_E)\in L_i-\{b\}$. Similarly we have an exit 
vertex $y$ of an endblock in $L_i$ with $p(y)\in K-\{a\}$. Now by Proposition \ref{int} we 
have $x_E\neq p(y)$ and both are contained in $S_i$, proving the claim. 

Assume that $L_1,\ldots ,L_t$ are labelled so that for $i\in [1,t']$, $|S_i|\geq 2$ 
and $|S_i|=1$ for $i\in [t'+1,t]$. By the previous paragraph we have $t'\geq 1$. For 
all $I\subset [t]$ we let $S_I=\bigcup _{i\in I} S_i$. 

Beginning with the $\{S_1,\ldots ,S_{t'}\}$, repeatedly replace sets $S_I$ and $S_J$ in 
this list with $S_{I\cup J}$ if $|\mbox{span}_{G_a}(S_I)\cap \mbox{span}_{G_a}(S_J)|\geq 2$.
When this proceedure ends we are left with sets $\{S_{I_1},\ldots ,S_{I_p}\}$. 

Now clearly $\mbox {span}_{G_a}(S_{I_l})$ is a union of blocks of $K$ for all $l\in [p]$ and 
by our construction proceedure above, no two can share a block. Also from the claim above, each 
endblock $E$ of $K$ is contained in $\mbox {span}_{G_a}(S_{I_l})$ for some $l\in [p]$. Combining 
these two facts it is easy to see that there is some $l\in [p]$ for which $\mbox{span}_{G_a}(S_{I_l})$ 
is separated from $G_a - \mbox{span}_{G_a}(S_{I_l})$ in $G_a$ by a single vertex $a'$. 

We are now ready to choose $J$. Let $N\subset [t'+1,t]$ consist of all $n$ for 
which $S_n\cap (\mbox{span}_{G_a}(S_{I_l}) -\{a'\})\neq \emptyset$ and let $M=I_l\cup N$. 
We take $J=G[\mbox{span}_{G_a}(S_M)\cup \mbox{span}_{G_b}(T_M)]$. 

Lastly we choose $b'$. We pick this vertex depending on whether $|M|=1$ or $|M|\geq 2$. 
If $|M|=1$ then $I_l=\{i\}$ for some $i\in [t']$ and $N=\emptyset$. Since 
$G_b[\mbox{span} _{G_b}(T_i)]$ is a connected union of blocks of $L_i$ 
and contains a vertex in the interior of every endblock of $L_i$, this graph must be 
separated from $G_b-G_b[\mbox{span} _{G_b}(T_i)]$ by a single vertex $b'$ in $G_b$. 

If $|M|\geq 2$ then as $C$ contains at least two limbs of $b$, by Lemma \ref{body} 
$\mbox{Body}(b)=\{b\}$. Then $J=G[\mbox{span}_{G_a}(S_{M})\cup (\bigcup _{j\in M} V(L_j))]$ 
and we may take $b'=b$. 

It is clear from construction that $J$ satisfies (i). It is also easy to show 
that $J$ is $2$-connected. Indeed, suppose we remove $v\in J\cap G_a$. 
Given any vertex $w\in J\cap G_a-\{v\}$, $w\in G[\mbox{span}_{G_a}(S_{I_l})]$ so either 
$w\in S_{I_l}$ or $w$ lies on a path between two elements of $S_{I_l}$. In both cases there 
exists a path from $w$ to $J\cap G_b$. Since $J\cap G_b$ is connected this shows that $J-v$ is 
connected for all $v\in J\cap G_a$. A similar argument shows that $J-v$ is connected for all
$v\in J\cap G_b$.

We now show that (ii) holds for $J$. Suppose that $v\in J-\{a',b'\}$ and that we are looking for 
an $a'-v$ path not containing $b$ of length at least $2^{d-2}$. We claim the following:
\vspace{2mm}

\noindent \textbf{Claim 1:} $J-b$ contains a $2$-connected subgraph $J'$ containing all of $J\cap G_a$.
\vspace{2mm}

If $b\notin J$ then this is immediate taking $J'=J$, so we may assume that $b=b'\in J$. It suffices 
to show that there exists such a subgraph $J'$ of 
$G[\mbox {span}_{G_a}(S_{I_l})\cup \mbox{span}_{G_b}(T_{I_l})]$.

For each $i\in [t']$, $G[\mbox {span}_{G_a}(S_{i})\cup \mbox{span}_{G_b}(T_{i})]-\{b\}$ has a 
$2$-connected subgraph which contains all of $\mbox{span}_{G_a}(S_i)$, namely 
$G[\mbox {span}_{G_a}(S_{i})\cup \mbox{span}_{G_b}(p(S_i))]-\{b\}$ where $p(S_i):=\{p(s): s\in S_i\}$. 
Here $|S_i|\geq 2$ guarantees that this graph is $2$-connected.
Now since the union of two $2$-connected graphs which intersect at least two points is still 
a $2$-connected graph, from the joining proceedure which produced the set $I_l$ we must have 
that the union of $G[\mbox {span}_{G_a}(S_{i})\cup \mbox{span}_{G_b}(p(S_i))]-\{b\}$ 
for $i\in I_l$ is a $2$-connected graph $J'$. Moreover, this graph clearly contains all of 
$G[\mbox{span}_{G_a}(S_{I_l})] = J\cap G_a$. This proves the claim.

We now use $J'$ to find the $a'-v$ path claimed in (ii). First find a path $P_1$ from $v$ to 
some $w\in J'-a'$ which avoids $b$. Such a path is immediate if $v\in J'$, so we may assume 
$v\notin J'$. Let $v\in L_i$, $i\in M$. If $L_i-b$ contains some element of $J'$, take 
$P_1$ to be the shortest path in $L_i-b$ from $v$ to an element $w$ in $J'$. If not, we take 
$P_1'$ to be a path in $L_i-b$ to an exit vertex $x_F$ of some endblock $F$ of $L_i$, 
$p(x_F)=w$ and $P_1=P_1'x_Fp(x_F)$. Note that in both these cases $P_1$ intersects $J'$ only 
in one vertex $w\neq a'$.

Now take an endblock $E$ of $J\cap G_a$. If $w\in \mbox{int}(E)$ then $G_a$ contains a path 
from $w$ to $a'$ which extends an endblock path in $E$. As such a path has length at least $2^{d-2}$ 
we can assume $w\notin \mbox{int}(E)$. Let $J''$ denote the graph $J'$ with $\mbox{int}(E)$ 
contracted to a single vertex $e$. It is easy to see that this graph is still $2$-connected. 
Therefore there exists two vertex disjoint paths $P_2$ and $P_3$ from $\{a',w\}$ to 
$\{\mbox{cutv}(E),e\}$. Say that these paths are $P_2$ from $a'$ to $\mbox{cutv}(E)$ and $P_3$ 
from $w$ to $e$. If $P_3=P_3'xe$, $x$ must have a partner $p(x)\in \mbox{int}(E)$. This gives 
a path $P_3'xp(x)$ from $w$ to $p(x)$ in $J$. Now since $E$ is a $2$-connected graph and for all 
$v\in E-\{\mbox{cutv}(E),p(x)\}$ $d_E(v)\geq d-1$, we can apply Theorem \ref{connmod} to $E$ to 
find a $p(x)-\mbox{cutv}(E)$ path $P_4$ of length at least $2^{d-2}$. Combining all of these paths 
gives a $v-a'$ path $P=P_1P_3'xp(x)P_4P_2^{(r)}$ of length at least $2^{d-2}$, where $P_2^{(r)}$ 
is $P_2$ reversed. 

An identical argument gives the $b'-v$ path claimed in (ii).
}
\end{proof}

\begin{figure}
\centering
\input{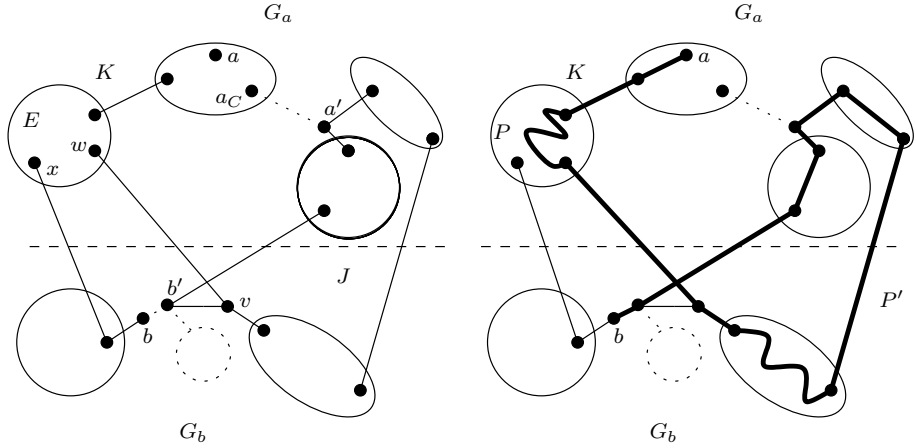}
\caption{An illustration of the case $w\in \mbox{int}(E)$ in the proof of
Lemma \ref{splittingdirection}. The broken dotted pieces represent vertices in $G_C$
that are left out of $J$.}
\label{together}
\end{figure}

We are now ready to prove Lemma \ref{splittingdirection}.

\vspace{2mm}

\noindent \emph{Proof of Lemma \ref{splittingdirection}.}
Suppose for contradiction that $G$ does not contain an $a-b$ path of length at least $2^{d-1}$. 
Then by Corollary \ref{limbcomponent}, the interaction digraph $H$ of $G$ must contain a connected 
component $C$ consisting entirely of limbs.

As $G$ does not contain an $a-b$ path of length at least $2^{d-1}$, we can apply Lemma 
\ref{nice} to find a $2$-connected subgraph $J$ of $G_C$ and vertices $a'$ and 
$b'$ which satisfy Lemma \ref{nice} (i) and (ii). Note that $|J|<|G|$ since $H$ 
contains two connected components and $J$ is contained entirely in one of them.

Now if there are no edges between $J-\{a',b'\}$ and $G-J$, all $v\in J-\{a',b'\}$ have degree at 
least $d$ in $J$. But then since $|J|<|G|$, $J$ contains an $a'-b'$ path $P$ of length at least 
$2^{d-1}$. Extending this path from $a'$ to $a$ and from $b'$ to $b$ gives an $a-b$ path 
of length at least $2^{d-1}$, a contradiction. Therefore such an edge must exist, joining say 
$v\in J-\{a',b'\}$ to $w\in G-J$. By Lemma \ref{nice} (i) $w\notin G_C$. We may assume $w\in G_a$. 

Suppose first that $w\in \mbox{Body}(a)$. Here $\mbox{Body}(a)\neq \{a\}$ by 
Lemma \ref{nice} (i). Take an $a'-v$ path of length $2^{d-2}$ in $J$ as guaranteed by Lemma 
\ref{nice} (ii) which does not contain $b$. This path extends in $G_C$ to an $a_C-v$ path $P_1$, 
where $a_C=G_C\cap \mbox{Body}(a)$. As this path lies entirely in $G_C$, it can only intersect 
$\mbox{Body}(a)$ in $a_C$ and $\mbox{Body}(b)$ in at most $b_C=G_C\cap \mbox{Body}(b)$. Pick a 
limb $K$ of $a$ not contained in $C$ and an endblock $E$ of $K$. $K$ contains a path of length 
$2^{d-2}$ from $\mbox{Joint}(K)$ to the exit vertex $x_E$ of $E$. As $p(x_E)\notin G_C$, we can 
find a path $P_2$ from $p(x_E)$ to $b$ in $G_b$ disjoint from $P_1$. But now since $\mbox{Body}(a)$ 
is $2$-connected we can find vertex disjoint paths from $\{a, \mbox{Joint}(K)\}$ to $\{a_C, w\}$. 
Combining these paths with $P_1$ and $P_2$ we obtain a path of length at least $2^{d-1}+2$ from $a$ 
to $b$, again a contradiction. 

Therefore we can assume $w\in K$ for some limb $K$ of $a$ not in $C$. If $w\notin \mbox{int}(E)$ 
for some endblock $E$ of $K$ then we can proceed exactly as in the case $w\in \mbox{Body}(a)$ 
above to find an $a-b$ path of length at least $2^{d-1}$, so we may assume $w\in \mbox{int}(E)$. 
Then $K$ contains a path of length at least $2^{d-2}$ from $w$ to $\mbox{Joint}(K)$. Joining 
this path to the $v-b'$ path guaranteed by Lemma \ref{nice} (ii) via the edge $wv$ we obtain a 
$\mbox{Joint}(K)-b'$ path of length at least $2^{d-1}+1$. Extending this path from $\mbox{Joint}(K)$ 
to $a$ and from $b'$ to $b$ we again find an $a-b$ path of length at least $2^{d-1}+1$. This 
contradicts our assumption and proves the Lemma. \hspace{10.3cm} $\square $

% Section 6: Removing the Degree Assumption

\section{Removing the Degree Assumption}

Again, let $G$ be a $2$-connected subgraph of $Q_n$ with $a,b\in G$ such that 
$d_G(v)\geq d$ for all $v\in V(G)-\{a,b\}$. Also, suppose that Theorem 
\ref{connmod} holds for smaller degrees and for all graphs $G'$ with 
$|G'|<|G|$.

If we could find a splitting direction $i$ for $G_a$ and $G_b$ in Lemma \ref{splitting} 
so that $d_{G_a}(a)\geq 2$ and $d_{G_b}(b)\geq 2$ then using Lemma \ref{splittingdirection} 
$G$ would contain an $a-b$ path of desired length. This is possible if $d_G(a)\geq 3$ and 
$d_G(b)\geq 3$, so we may assume that say $d_G(a)=2$. However, we can not guarantee 
this in general -- for example, $a$ and $b$ could be adjacent with both having only one 
other neighbour in $G$.

Now the condition $d_{G_a}(a)\geq 2$ and $d_{G_b}(b)\geq 2$ in previous sections 
ensured that all endblocks of $G_a$ and $G_b$ contained long paths, which is false 
if say $a$ has a single neighbour $a'$ in $G_a$. This fact was then used in 
Lemma \ref{list} (iv) to show that $a$ has at least two limbs in $G_a$ which in turn 
was crucially used numerous times in our analysis of $H$ e.g. Lemma \ref{body}. 

In this section, we will extend the arguments from the proof of Lemma \ref{splittingdirection} 
to prove Theorem \ref{connmod}. The first step is the following lemma.

% Section 6: First Lemma

\begin{lem} 
\label{remove}
Let $G$ be a $2$-connected subgraph of $Q_n$ with $a,b\in G$  
such that $d_G(a)=2$ and $d_G(v)\geq d$ for all $v\in V(G)-\{a,b\}$, 
where $d\geq 3$.
Suppose that Theorem \ref{connmod} holds for smaller degrees and 
for all graphs $G'$ with $|G'|<|G|$. 
Furthermore, suppose that $G$ does not contain an $a-b$ path 
of length at least $2^{d-1}$. 
Then the following hold:
\begin{itemize}
{
\item[\emph{(i)}] $G-a$ is a $2$-connected graph

\item[\emph{(ii)}] There is a splitting direction $i$ for $G_a$ and $G_b$ 
so that $d_{G_a}(a)\geq 2$ or $d_{G_b}(b)\geq 2$.
}
\end{itemize}
\end{lem}

% Section 6: Proof of First Lemma

\begin{proof}
 {

(i) If $G-a$ is not $2$-connected, it has at least two endblocks in its
block-cutvertex decomposition, one of which $E$ has $b\notin \mbox{int}(E)$. Now
since $G$ is $2$-connected, $a$ must be joined to the interior of all the
endblocks of $G-a$. As $d_G(a)=2$, $G-a$ has exactly two endblocks with $a$
having exactly one neighbour in the interior of each. Let $w$ be this neighbour
in $E$.

Now $E$ is $2$-connected (as $d\geq 3$) and all vertices in $E-\lbrace
\mbox{cutv}(E), w\rbrace $ have degree at least $d$ in $E$. Since $|E|<|G|$,
Theorem \ref{connmod} applies to give a path $P$ of length at least $2^{d-1}$ from $w$ 
to $\mbox{cutv}(E)$. Extending this path on either side to $a$ and $b$ 
respectively, we have an $a-b$ path of length at least $2^{d-1}$, a 
contradiction.

(ii) We can always choose such a direction if $a$ and $b$ are at Hamming
distance at least three in $Q_n$ or if one of $a$ or $b$ have degree greater
than 2 in $G$. So we can assume $a$ and $b$ are at Hamming distance one or 
two and both have degree exactly two in $G$. 

First consider $a$ and $b$ at Hamming distance one. If they are not adjacent in
$G$ we can choose the direction on which they differ for $i$ so we can assume
they are adjacent. Then $a$ and $b$ both have one other neighbour in $G$, $a'$
and $b'$ respectively. Now if $G-\{a,b\}$ is $2$-connected we can apply Theorem
\ref{connmod} to $G-\{a,b\}$ with $a'$ and $b'$ in place of $a$ and $b$. This 
gives an $a'-b'$ path of length at least $2^{d-1}$. Adjoining the edges $aa'$ and
$bb'$ to this path we have an $a-b$ path of length $2^{d-1}+2$, more than
enough. If $G-\{a,b\}$ is not $2$-connected it is easily seen that $a'$ and $b'$
must lie in the interior of different endblocks of $G-\{a,b\}$. We can therefore
find a path from $a'$ to $b'$ in $G-\{a,b\}$ which extends two endblock paths.
Adjoining the edges $aa'$ and $bb'$ to this path, we have an $a-b$ path of
length at least $2^{d}+2$, a contradiction.

If $a$ and $b$ are at Hamming distance two, we can always find such a direction 
$i$ unless $a$ and $b$ are joined to the same two neighbours in $G$,
$a'$ and $b'$ say. Then $\lbrace a,a',b,b'\rbrace $ form a $C_4$ with $a$
opposite $b$. Working with $G-\lbrace a,b\rbrace $, $a'$ and $b'$ as above, we
again obtain an $a-b$ path of desired length in $G$.
}
\end{proof}

We are now ready for the proof of Theorem \ref{connmod}.
\vspace{2mm}

% Proof of Theorem \ref{connmod} 2^{d-1}

\noindent \emph{Proof of Theorem \ref{connmod}.} 
The proof is by induction on $d$ and $|G|$. The base case $d=2$ follows from Menger's theorem, 
as if $G$ is $2$-connected it contains two disjoint $a-b$ paths, one of which must have length 
at least $2$. 

Suppose that $G$ is as in the statement of the Theorem and that Theorem \ref{connmod} holds for 
all smaller values of $d$ and for all graphs $G'$ with $|G'|<|G|$. Suppose for contradiction that 
$G$ does not contain a $a-b$ path of length at least $2^{d-1}$. By Lemma \ref{splittingdirection}, 
we must have that for all choices of $i$ in Lemma \ref{splitting} either $d_{G_a}(a)\leq 1$ or 
$d_{G_b}(b)\leq 1$. 

Take the splitting direction $i$ for $G_a$ and $G_b$ as in Lemma \ref{remove}(ii). We will 
assume without loss of generality that $d_{G_b}(b)\geq 2$. By the previous paragraph we must 
have $d_{G_a}(a)=1$. Let the neighbours of $a$ be $a'\in G_a$ and $v\in G_b$. 

Now we can assume $v\neq b$. Otherwise, by Lemma \ref{remove}(i), $G-a$ is $2$-connected and as all 
$v\in V(G-a)-\{a',b\}$ have degree at least $d$ in $G-a$, by induction, $G-a$ contains an $a'-b$ 
path of length at least $2^{d-1}$. Appending the edge $aa'$ to this path we obtain an $a-b$ path 
of length at least $2^{d-1}$, a contradiction.

% Claim: v\in int(E_v)

Lemma \ref{list} (i)-(iv) still hold for $G_b$ with the same proofs as before. In particular, 
$b$ still has at least two limbs. We make the following claim:
\vspace{2mm}

\noindent \textbf{Claim 2:} $v\in \mbox{int}(E_v)$ for some endblock $E_v$ of $G_b$
\vspace{2mm}

Suppose otherwise. From Lemma \ref{list}(iii) $b$ does not lie in the interior of an endblock of 
$G_b$ and by Lemma \ref{list}(iv) $G_b$ contains two vertex disjoint paths $P_1$ from $v$ to 
$\mbox{cutv}(E_1)$ and $P_5$ from $\mbox{cutv}(E_2)$ to $b$, where $E_1$ and $E_2$ are two endblocks 
of $G_b$. Taking exit vertices $x_1$ and $x_2$ of $E_1$ and $E_2$ respectively, by induction on $d$, 
$E_1$ contains a path $P_2$ of length at least $2^{d-2}$ from $\mbox{cutv}(E_1)$ to $x_1$ and 
$E_2$ contains a path $P_4$ of length at least $2^{d-2}$ from $x_2$ to $\mbox{cutv}(E_2)$. 
Taking a path $P_3$ from $p(x_1)$ to $p(x_2)$ in $G_a-a$ and combining the paths, $G$ contains an 
$a-b$ path $avP_1P_2x_1p(x_1)P_3p(x_2)x_2P_4P_5$ of length at least $2^{d-1}$. This contradicts 
our assumption and proves the claim.

We now again construct an interaction digraph $H$ but this time it is built from the
limbs of $a'$ and $b$ instead of those of $a$ and $b$. Note that $\{a,a'\}$ is a
limb of $a'$ and so, both $a'$ and $b$ have at least two limbs. Take
$H=(A',B,\overrightarrow{E})$ to be a bipartite multidigraph on vertex sets
$A'=\{K_1,\ldots ,K_r\}$ and $B=\{L_1,\ldots L_s\}$, the set of limbs of $a'$
and $b$ respectively. We also adjoin $\mbox{Core}(b)$ to $B$ if it is non-empty
($\mbox{Core}(a')=\emptyset$ since $a'$ is a cutvertex of $G_a$). Now each
endblock of $G_a$ other than $\{a,a'\}$ contains at least two exit vertices,
as in Lemma \ref{list}(i). Therefore for each endblock $E$ of $G_a$ or $G_b$
other than $\{a,a'\}$ we can pick an exit vertex $x_E$ with $p(x_E)\neq a',b$.
From our claim above we can pick $x_{E_v}=v$. 
Now adjoin a directed edge from $K\in A'$ to $L\in B$ for each endblock $E$ in 
$K$ with $p(x_E)\in L$ and a directed edge from $L\in B$ to $K\in A'$ for each 
endblock $E$ in $L$ with $p(x_E)\in K$. Note that every limb other than $\{a,a'\}$ 
still has an outneighbour in $H$. 

For this $H$ Lemma \ref{paths3} and Lemma \ref{body} still hold, again with the same
proofs as before. Using these two, as in Corollary \ref{limbcomponent}, we can show 
that $H$ contains a connected component $C$ consisting entirely of limbs which does 
not contain the limb $\{a,a'\}$. 
Indeed, $b$ has at least two limbs so pick one, $L\in B$, not containing $E_v$ 
and take $C$ to be the connected component of $H$ containing $L$. As $v$ is the 
unique neighbour of $a$ in $G_b$ and $v\notin L$, if $\{a,a'\}\in C$ then $H$ would
contain a path of length three, contradicting Lemma \ref{paths3}. Furthermore,
since $\mbox{Core}(a')=\emptyset $, if $C$ did not consist entirely of limbs of
$a'$ and $b$, $\mbox{Body}(b)\neq \{b\}$ and $C$ contains two vertices of $B$,
contradicting Lemma \ref{body}. 

The remainder of the proof of Theorem \ref{connmod} is almost identical to that of Lemma 
\ref{splittingdirection}. We can apply Lemma \ref{nice} to find a subgraph $J$
of $G_C$. Using this subgraph as in the proof of Theorem \ref{splittingdirection} 
we either obtain an $a'-b$ path of length at least $2^{d-1}$ which is contained 
entirely in $G_C$ or an $a'-b$ path of length at least $2^{d-1}+1$ in $G$. In the 
first case we find our $a-b$ path by appending the edge $a'a$ the to $a'-b$ path. 
In the second case, unless $a$ is already a vertex of this path we can also do this. 
But if $a$ a vertex of this $a'-b$ path, it must occur as the second vertex. Deleting 
$a'$ from the path, we obtain an $a-b$ path of length at least $2^{d-1}$, as required.
\hspace{6.5cm} $\square$

% Section 7: A Tight Bound

\section {A Tight Bound}

In this section we will prove Theorem \ref{conn}. Its proof has the same structure
as Theorem \ref{connmod} but requires more care in various arguments. The proof will 
again be by induction on $d$. 
The base case $d=2$ is immediate unless $a$ and $b$ are at Hamming distance $2$ apart. 
If this is the case and $G$ is not isomorphic to $Q_2$ pick any vertex $v$ of $G$ not 
in the unique $2$-cube containing $a$ and $b$. By $2$-connectivity $G$ contains vertex 
disjoint $a-v$ and $v-b$ paths, which when combined give a path of length at least $3$, 
as required.

We will suppose that the Theorem fails for some $d>2$ and take $G$ to be a minimal 
counterexample so that Theorem \ref{conn} holds for all smaller degrees and all 
graphs $G'$ with $|G'|<|G|$. To begin we will prove the analogue of Lemma 
\ref{splittingdirection}.

% Section 7: Analogue of Key Lemma - degGa(a)\geq 2, degGb(b)\geq 2.

\begin{lem}
\label{splittingdirectionextension}
Let $G$ be a $2$-connected subgraph of $Q_n$, not isomorphic to $Q_d$ and let $a,b\in V(G)$ such 
that $d(v)\geq d$ for all $v\in V(G)-\{a,b \}$, where $d\geq 3$. Suppose that Theorem \ref{conn} 
is true for smaller degrees and all graphs $G'$ with $|G'|<|G|$. Suppose furthermore that there 
exists a splitting direction $i$ for $G_a$ and $G_b$ in Lemma \ref{splitting} for which 
$d_{G_a}(a)\geq 2$ and $d_{G_b}(b)\geq 2$. Then $G$ contains an $a-b$ path of length at least 
$2^{d}-1$.
\end{lem}

% Section 7: First Lemma List

Our first step in the proof of Lemma \ref{splittingdirectionextension} is to establish the 
analogue of Lemma \ref{list}.

\begin{lem}
{
\label{listextension}
Let $G, a, b, G_a$ and $G_b$ be as in the statement of Lemma \ref{splittingdirectionextension}. 
If $G$ does not contain an $a-b$ path of length at least $2^d-1$ the following hold:
\begin{itemize}
{
\item[\emph{(i)}] Every endblock $E$ of $G_a$ which does not contain $a$ in its
interior contains at least two exit vertices $x$ and $x'$. Furthermore, we can 
choose these so that $E$ contains $\mbox{cutv}(E)-x$ and $\mbox{cutv}(E)-x'$ 
paths of length at least $2^{d-1}-1$. 

\item[\emph{(ii)}] $G_a$ is not $2$-connected.

\item[\emph{(iii)}] $a$ does not lie in the interior of an endblock in $G_a$.

\item[\emph{(iv)}] $a$ must have at least two limbs.
}
\end{itemize}
}
\end{lem}

% Section 7: Proof of Lemma List Extension

\begin{proof}
 {

(i) Here the proof of Lemma \ref{list}(i) needs only a small change. If $E$ is isomorphic 
to $Q_{d-1}$ we can choose any two neighbours of $\mbox{cutv}(E)$ for $x$ and $x'$, 
so we may assume $E$ is not isomorphic to $Q_{d-1}$. Now $G$ is $2$-connected so $E$ contains 
at least one exit vertex $x$. Suppose for contradiction that this is the only one. 
Then as $E$ is $2$-connected, not isomorphic to $Q_{d-1}$ with $d_E(v)\geq d$ for all 
$v\in E-\{\mbox{cutv}(E), x\}$ and $|E|<|G|$, it contains a $\mbox{cutv}(E)-x$ path of length 
at least $2^d-1$. Extending this path as before we obtain an $a-b$ path of length at least 
$2^d-1$, a contradiction. Therefore $E$ has two exit vertices and as $E$ is not isomorphic to 
$Q_{d-1}$, $E$ contains paths of length at least $2^{d-1}-1$ from $\mbox{cutv}(E)$ to both of 
them.

(ii) The change to the proof of the Lemma \ref{list} (ii) in this case is a little 
more involved. Suppose for contradiction that $G_a$ is $2$-connected.

First suppose that $G_b$ is not $2$-connected. If there exists an endblock $E$ in
$G_b$ such that $b\notin E$, we have a path $P_1$ in $G_b$ of length at least $1$ from
$b$ to cutv($E$). Pick an exit vertex $x$ of $E$ such that $E$ contains a $\mbox{cutv}(E)-x$ 
path $P_2$ of length at least $2^{d-1}-1$ and $p(x)\neq a$ -- this exists by (i). 
Combining the paths $P_1$ and $P_2$ above with the $a-p(x)$ path of length
$2^{d-1}-2$ in $G_a$ guaranteed by Theorem \ref{conn}, $G$ contains
an $a-b$ path of length at least $1+(2^{d-1}-1)+1+(2^{d-1}-2)=2^d-1$, a
contradiction. So if $G_b$ is not $2$-connected $b$ must lie in \emph{every} 
endblock $E_1,\ldots ,E_t$ of $G_b$. Note that since $t\geq 2$ this implies 
$b\notin \mbox{int}(E_i)$ for any $i$. 

Now using (i) as with $E$ above, $E_1$ must have an exit vertex $x$ such that $E_1$ 
contains a $x-b$ path of length at least $2^{d-1}-1$, with $p(x)\neq a$. If $G_a$ 
were not isomorphic to $Q_{d-1}$, it contains a path of length $2^{d-1}-1$ from $a$ to
$p(x)$. Combining these two with the edge $xp(x)$ we obtain an $a-b$ path of
length $2^{d}-1$. 

Therefore we can assume $G_a$ is isomorphic to $Q_{d-1}$. Then $G_a$ contains a 
path $P_3$ of length at least $2^{d-1}-1$ from $a$ to any of its neighbours. Take 
a neigbhour $x$ such that $p(x)\neq b$. Here $p(x)$ must be in $\mbox{int}(E_i)$ 
for some $i\in [t]$. Now $t\geq 2$ so $E_i$ is not isomorphic to $Q_{d-1}$ 
-- otherwise $G_a$ would receive too many edges from $G_b$ by (i) above. Since 
Theorem \ref{conn} holds for smaller degrees, $E_i$ contains a path $P_4$ from 
$b=\mbox{cutv}(E_i)$ to $p(x)$ of length at least $2^{d-1}-1$. Combining $P_3$ 
and $P_4$ with the edge $xp(x)$ we have an $a-b$ path of length at least $2^d-1$, 
a contradiction.

The case when $G_b$ is $2$-connected is very similar. We can obtain two paths of
length at least $2^{d-1}-1$ in $G_a$ and $G_b=E_1$ if neither of the two are
isomorphic to $Q_{d-1}$ and if one is isomorphic to $Q_{d-1}$ we can use the same
argument as in the case where $G_a$ is isomorphic to $Q_{d-1}$ and $t\geq 2$
above.

(iii) This is similar to (ii) but a little easier. Suppose that $a$ is contained 
in the interior of some endblock $E$ of $G_a$. As Theorem \ref{conn} holds 
for degrees smaller than $d$, $E$ contains an $a-\mbox{cutv}(E)$ path $P_1$ of 
length at least $2^{d-1}-2$. Since $G_a$ is not $2$-connected by (ii), it also 
contains a second endblock $E'$. Now $E'$ must contain an exit vertex $x$ with 
$p(x)\neq b$ for which $E'$ contains a $\mbox{cutv}(E')-x$ path $P_3$ of length 
at least $2^{d-1}-1$. Joining $\mbox{cutv}(E)$ to $\mbox{cutv}(E')$ by a path 
$P_2$ in $G_a$ and $p(x)$ to $b$ with a path $P_4$ in $G_b$ gives an $a-b$ path 
$P=P_1P_2P_3xp(x)P_4$ of length at least $(2^{d-1}-2)+0+(2^{d-1}-1)+1+1=2^d-1$, 
a contradiction.

(iv) Again follows from (ii) and (iii) as in Lemma \ref{list}(iv).
}
\end{proof}

% Section 7: Discussion between Lemma List Ext and Lemma Int Ext

The above modifications demonstrate the main problem in moving from the bounds of
Theorem \ref{connmod} to bounds of Theorem \ref{conn} -- on combining endblock paths together 
without any care as before, we are usually left short a small number of vertices. 
While in the above Lemma we were able to exploit some small extremal arguments to 
obtain these extra vertices, such arguments do not allow us to prove the natural 
analogue of Proposition \ref{int}. Indeed, we may have an endblock $E$ of $G_a$ 
not isomorphic to $Q_{d-1}$ and an endblock $F$ of $G_b$ isomorphic to $Q_{d-1}$
with $\mbox{cutv}(E)=a$ and $\mbox{cutv}(F)=b$. If there is a vertex $x\in \mbox{int}(E)$ 
at odd Hamming distance from $a$ adjacent to a vertex $y\in \mbox{int}(F)$ at even 
Hamming distance from $b$, our induction hypothesis only allows us to find a path of 
length at least $(2^{d-1}-1)+1+(2^{d-1}-2)=2^d-2$ from $a$ to $b$.

Now Lemma \ref{int} was important in the proof Theorem \ref{connmod}. In particular, 
it was used in the construction of the interaction digraph of $G$ and crucially in 
the proof of Lemma \ref{nice} where it guaranteed the existence of the $2$-connected 
subgraph $J$. The next proposition is a weakened version of Proposition \ref{int} 
which will play a similar role in the proof of Theorem \ref{conn}.

From Lemma \ref{listextension}(i) above, for every endblock $E$ from $G_a$ we can 
pick an exit vertex $x_E$ so that $E$ contains a $\mbox{cutv}(E)-x_E$ path of length 
at least $2^{d-1}-1$ and for which $p(x_E)\notin \{a,b\}$. Similarly pick such 
exit vertices $y_F$ for endblocks $F$ in $G_b$.

% Section 7: Proposition Int Ext

\begin{prop}
\label{intextension}
Let $G, a, b, G_a$ and $G_b$ be as in the statement of Lemma \ref{splittingdirectionextension}. 
Suppose $G$ does not contain an $a-b$ path of length at least $2^d-1$. Then for every 
endblock $E$ of $G_a$, $p(x_E)\neq y_F$ for any endblock $F$ of $G_b$. Furthermore, if 
$\mbox{Body}(a)\neq \{a\}$ or $\mbox{Body}(b)\neq \{b\}$ we also have 
$p(x_E)\notin \mbox{int}(F)$.
\end{prop}

% Section 7: Proof of Proposition Int Ext

\begin{proof}
From Lemma \ref{listextension}(iii) $a\notin \mbox{int}(E)$ and $b\notin \mbox{int}(F)$. Suppose 
for contradiction that $p(x_E)=y_F$. Combining the $\mbox{cutv}(E)-x_E$ path in $E$ with 
the $\mbox{cutv}(F)-y_F$ path in $F$ guaranteed by Lemma \ref{listextension}(i) via the 
edge $x_Ep(x_E)=x_Ey_F$, we have a $\mbox{cutv}(E)-\mbox{cutv}(F)$ path of length at least 
$(2^{d-1}-1)+1+(2^{d-1}-1)=2^d-1$. As this path extends to an $a-b$ path, we have a 
contradiction.

The second part is similar. Since $F$ is $2$-connected and we have $d_F(v)\geq d-1$ for all 
$v\in F-\{\mbox{cutv}(F),p(x_E)\}$, $F$ contains a $p(x_E)-\mbox{cutv}(F)$ path of length 
at least $2^{d-1}-2$. Combining this with the $\mbox{cutv}(E)-x_E$ path of length at least 
$2^{d-1}-1$ in $E$ via the edge $x_Ep(x_E)$ we have a $\mbox{cutv}(E)-\mbox{cutv}(F)$ path 
of length at least $2^d-2$. If $\mbox{Body}(a)\neq \{a\}$ or $\mbox{Body}(b)\neq \{b\}$, 
extending this path to an $a-b$ path takes at least one more edge, again giving an $a-b$ 
path of length at least $2^d-1$, a contradiction.
\end{proof}

% Section 7: Construction of H

We now look towards an slightly altered construction for $H$. Let our interaction 
digraph $H= \{A,B,\overrightarrow E\} $ again be a bipartite multidigraph whose 
bipartition consists of the limbs of $a$ and $b$ respectively. Again we additionally 
adjoin Core$(a)$ and Core$(b)$ to $A$ and $B$ respectively if they are non-empty. 
For each endblock $E$ of $G_a$, adjoin a directed edge to $H$ from $K\in A$ to $L\in B$ if 
$E$ is an endblock of $K$ and $p(x_E)\in L$. Similarly, for each endblock $F$ of $G_b$, 
adjoin a directed edge to $H$ from $L\in B$ to $K\in A$ if $F$ is an endblock of $L$ and 
$p(y_F)\in K$. Note again that every limb in $H$ has outdegree at least $1$ as it contains 
an endblock.

% Section 7: Lemma Paths3 Extension

We now prove the analogue of Lemma \ref{paths3}. The original proof is complicated 
by the fact that endblock paths can join into the interior of endblocks, as 
discussed above.

\begin{lem} 
{
\label{path3extension}
Let $G, a, b, G_a$ and $G_b$ be as in Lemma \ref{splittingdirectionextension}. Suppose that 
$G$ does not contain an $a-b$ path of length at least $2^d-1$. Then $H$ cannot contain an 
undirected path of length three.
}
\end{lem}

% Proof of Lemma Paths3 Extension

\begin{proof}
{
The proof of Lemma \ref{paths3} applies unchanged if we can guarantee that for any endblock 
$E$ of $G_a$ or $G_b$, $p(x_E)\notin \mbox{int}(F)$ for any endblock $F$ of $G_b$ or $G_a$. 
Using Proposition \ref{intextension} above we can therefore focus on the case where 
$\mbox{Body}(a)=\{a\}$ and $\mbox{Body}(b)=\{b\}$ i.e. $a$ is a cutvertex of $G_a$ and 
$b$ is a cutvertex of $G_b$.

Suppose for contradiction that $H$ contains a path of length at least $3$. 
If one of the interior vertices on this path has two out-neighbours in $Q$ 
the same argument as in the original proof will create a path between two exit 
vertices in this limb which extends two endblock paths of length $2^{d-1}-1$.
This gives a path of length at least $2^d-2$. Extending this path through $Q$ 
as in Theorem \ref{paths3} gives us an $a-b$ path of length at least $2^d-1$. 
Similarly, if both endvertices on this path have out-neighbours in $Q$ (that 
is $\overrightarrow{V_0V_1}$ and $\overleftarrow{V_2V_3}$ are edges of $Q$) 
we can find paths of length at least $2^{d-1}-1$ in both $V_0$ and $V_3$, 
which again can be joined through $Q$ to give an $a-b$ path of length at 
least $2^d-1$. This just leaves the case of a directed path

\begin{equation}
{
\label{description}
\overrightarrow{V_0V_1}, \overrightarrow{V_1V_2} \mbox{ and }
\overrightarrow{V_2V_3}.
}
\end{equation}

While here we obtain a path of length at least $2^{d-1}-1$ from the edge 
$\overrightarrow{V_0V_1}$ as before, in $V_1,V_2,V_3$ we might not be able to 
guarantee a full endblock path. Indeed, there is now the possibility
that the edges $\overrightarrow{V_{i-1}V_i}$ and $\overrightarrow{V_iV_{i+1}}$
correspond to an edge entering an endblock $E$ by a vertex $x$ in its interior
and the other edge leaving $E$ by a vertex $y$ in its interior. This does not
allow us to apply Theorem \ref{conn} to $E$ as $\mbox{cutv}(E)\in E-\{x,y\}$
may have degree lower than $d-1$ in $E$. If, however, this does not happen at one 
of $V_1$ or $V_2$ the same proof applies.

We may also assume that the exit vertex $x$ of $V_2$ guaranteed by $\overrightarrow{V_{2}V_3}$
has $p(x)\notin \mbox{int}(F)$ for any endblock $F$ in the interior of an endblock 
of $V_3$. Otherwise we can find an endblock path in $V_0$ of length at least $2^{d-1}-1$ 
and one in $V_3$ of length at least $2^{d-1}-2$. Since joining both of these through 
$V_1$ and $V_2$ joins at least four more vertices onto these paths, we can extend them 
to form an $a-b$ path of length at least $2^{d}+1$, more than required.

But now take any outneighour of $V_3$ in $H$. Combined with our path $Q$ above it is 
easily seen we can obtain a path $Q'$ of length three which is either (i) not of the 
form (\ref{description}), or (ii) contains $V_3$ as an interior vertex and allows for 
a full endblock path to be built through it. In both cases we are done.
}
\end{proof}

% Section 7: Lemma Body Extension

The following gives an analogue of Lemma \ref{body}. The proof from Lemma \ref{body} applies 
unchanged, on noticing that by Proposition \ref{intextension}, exit vertices of endblocks of 
$G_a$ and $G_b$ again cannot have partners in the interior endblocks of $G_b$ and $G_a$.

\begin{lem}
 {
\label{bodyextension}
Let $G, a, b, G_a$ and $G_b$ be as in the statement of Lemma \ref{splittingdirectionextension}. 
Suppose $G$ does not contain an $a-b$ path of length at least $2^d-1$. Furthermore, 
suppose that $\mbox{Body}(a)\neq \{a\}$. Then no connected component of $H$ contains two vertices 
of $A$.
}
\end{lem}

% Section 7: Limb Component Corollary

Combining Lemma \ref{path3extension} with Lemma \ref{bodyextension} as in Corollary 
\ref{limbcomponent}, we obtain the following:

\begin{cor} 
{
\label{limbcomponentextension}
Let $G, a, b, G_a$ and $G_b$ be as above. Suppose that $G$ does not 
contain an $a-b$ path of length $2^d-1$. Then the interaction digraph $H$ of 
$G$ has at least two connected components, one of which $C$ consists entirely 
of limbs.
}
\end{cor}

% Section 7: Nice Extension

\begin{lem}
 {
\label{niceextension}
Let $G, a, b, G_a$ and $G_b$ be as in the statement of Lemma \ref{splittingdirectionextension}.
Suppose $G$ does not contain an $a-b$ path of length at least $2^d-1$. Then taking $C$ 
as in Corollary \ref{limbcomponentextension}, $G_C$ has a $2$-connected subgraph $J$ 
containing two vertices $a'\in G_a$ and $b'\in G_b$ with the following properties:

% Three properties of subgraph J

\begin{itemize}
{
\item[\emph{(i)}] every vertex $v\in J - \lbrace a',b'\rbrace$ has degree at
least $d-1$ in $J$ and all the neighbours of $v$ in $G_C$ are contained in $J$

\item[\emph{(ii)}] for any vertex $v\in (J-\lbrace a',b'\rbrace )\cap G_a$, $J$ 
contains an $a'-v$ path not containing $b$ of length at least $2^{d-1}-2$. 
Furthermore, if $v$ has a neigbhour outside of $J$, $J$ contains a $b'-v$ path 
not containing $a$ of length at least $2^{d-1}$

\item[\emph{(iii)}] for any vertex $v\in (J-\lbrace a',b'\rbrace )\cap G_b$, $J$ 
contains a $b'-v$ path not containing $a$ of length at least $2^{d-1}-2$. 
Furthermore, if $v$ has a neighbour outside of $J$, $J$ contains an $a'-v$ path 
not containing $b$ of length at least $2^{d-1}$. 
}
\end{itemize}
}

\end{lem}

% Section 7: Proof of Lemma Nice Extension

\begin{proof} 
 {
The proof of this lemma is almost identical to that of Lemma \ref{nice} with Lemma 
\ref{path3extension} and \ref{bodyextension} taking the place of Lemma \ref{paths3} 
and \ref{body}. 

The only change to the argument of the proof is that in order to guarantee that we 
have $|S_i|\geq 2$ for some $i\in [t]$, we cannot now guarantee that $p(x_E)$ does 
not lie in the interior of any endblock of $G_b$. Instead, if $p(x_E)\in F$ for 
some endblock $F$ of $G_b$, by Proposition \ref{intextension} $\mbox{int}(F)$ contains 
an exit vertex $y_F$ with $p(y_F)\neq x_E$. This again shows that $|S_i|\geq 2$ for some 
$i\in [t]$ and therefore that $t'\geq 1$.

The bounds in (ii) and (iii) follow from our new bounds on the length of endblock paths. 
Indeed, suppose $v\in (J-\{a',b'\})\cap G_b$ say. Both the $a'-v$ and $b'-v$ paths in 
Lemma \ref{nice}(ii) contain entire endblock paths and therefore have length at least 
$2^{d-1}-2$. This gives the bound on the $b'-v$ path claimed. To obtain the $a'-v$ 
path, first note that as $v$ has a neigbhour $w$ outside $J$, by (i) it must be 
outside $G_C$. But then $w$ must be $v$'s partner, i.e. $w=p(v)$. Now in the 
last paragraph of the proof of Lemma \ref{nice}, the endblock $E$ is on the opposite 
side of $J$ from $v$. But the path constructed in Lemma \ref{nice} combines a path from $v$ 
to $E$ with an endblock path in $E$. The first of these has to have length at least $2$ 
as $p(v)\notin E$ and the second has length at least $2^{d-1}-2$. Combining these we 
obtain an $a'-v$ path of length at least $2^{d-1}$, as claimed.
}
\end{proof}

We will now give the proof of Lemma \ref{splittingdirectionextension}.
\vspace{3mm}

% Proof of Lemma 2^d-1 for d_{G_a}(a)\geq 2 and d_{G_b}(b)\geq 2

\noindent \emph{Proof of Lemma \ref{splittingdirectionextension}.}
Suppose for contradiction that $G$ does not contain an $a-b$ path of length at 
least $2^d-1$. Then by Corollary \ref{limbcomponentextension} the interaction 
digraph $H$ of $G$ contains a component $C$ consisting entirely of limbs.

Now as $G$ does not contain an $a-b$ path of length at least $2^d-1$, we can apply 
Lemma \ref{niceextension} to find a $2$-connected subgraph $J$ of $G_C$ and vertices 
$a'$ and $b'$ which satisfy Lemma \ref{niceextension} (i), (ii) and (iii). Again 
$|J|<|G|$.

% Part of the proof in which a-b path is built from edge entering J from outside is built

Suppose first that $v\in J-\{a',b'\}$ has a neighbour $w$ outside of $J$.
Without loss of generality take $v\in G_b$. Then $w\notin G_C$ by Theorem
\ref{niceextension}(i) and so $w\in K$ for some limb $K$ of $a$ or $w\in
\mbox{Core}(a)$. We will first deal with the case where $w\in K$.

% w \notin int(E)

If $w\notin \mbox{int}(E)$ for some endblock $E$ of $K$ take the $w-x_E$ path
$P_2$ in $K$ of length at least $2^{d-1}-1$ given by Lemma \ref{listextension}(i). 
Combining this with the path $P_1$ given from Lemma \ref{niceextension}(ii) 
in $J$ from $a'$ to $v$ of length $2^{d-1}-2$ and the edge $x_Ep(x_E)$, we have 
a path $P_1vwP_2x_Ep(x_E)$ from $a'$ to $p(x_E)$ of length at least 
$(2^{d-1}-2)+1+(2^{d-1}-1)+1= 2^d-1$. But this path extends to a path from 
$a$ to $b$, a contradiction.

% w \in int(E)

If $w\in \mbox{int}(E)$ for some endblock $E$ of $K$, we would like to combine
the $\mbox{Joint}(K)-w$ path guaranteed by induction on $E$ with the $v-b'$ path
in $J$ as given by Lemma \ref{niceextension}(ii) using the edge $wv$. This path 
extends to an $a-b$ path but may only have length $(2^{d-1}-2)+1+(2^{d-1}-2)=2^d-3$, 
too little for us. 

Instead, look at an outneighbour of $K$ in $H$. Let $K\in C'$ for some connected 
component $C'\neq C$ of $H$. If this is outneighbour is a limb, then $C'$ must
consist entirely of limbs, by Lemma \ref{bodyextension}. Therefore since all limbs 
have at least one outneighbour in $H$, $C'$ contains a path of the form 
$\overleftarrow{KW}$ or of the form $KVW$ where $\overleftarrow{VW}$ is
an edge of $H$. This allows us to build a path $P$ from $w$ to $\mbox{Joint}(W)$ 
of length at least $2^{d-1}+1$ in $C'$. Combining $P$ with an appropriate path 
from Lemma \ref{niceextension}(ii) via the edge $vw$ we obtain a path that extends 
to an $a-b$ path of length at least $(2^{d-1}-2)+1+(2^{d-1}+1)=2^d$, as required 
-- take this path to be the $a'-v$ path if $W\in B$ or the $b'-v$ path if $W\in A$.

If the outneighbour of $K$ in $H$ is instead $\mbox{Core}(b)$, again using Lemma
\ref{niceextension}(ii) we can find a $b'-v$ path in $J$ of length at least $2^{d-1}-2$
which extends through $C'$ to give a $y-b'$ path $P_2$ of length at least
$2^{d-1}$, where $y\in \mbox{Body}(b)$. Now $H$ must have a third connected
component $C''$ containing a limb of $b$ since $b$ has at least two limbs and
only one element of $B$ can lie in a component by Lemma \ref{bodyextension}. This
component gives an $a-z$ path $P_1$ of length at least $2^{d-1}-2$ where again
$z\in \mbox{Body}(b)$ and $P_1$ and $P_2$ are disjoint. As in the proof of Lemma
\ref{paths3} we can join $P_1$ and $P_2$ together in $\mbox{Core}(a)$ with a
small use of 2-connectivity to give an $a-b$ path of length at least $2^d-1$ as
required. This completes the case when $w\in K$. The case where $w\in \mbox{Core}(a)$ 
follows a similar argument, modifying the corresponding part of the proof of Lemma 
\ref{splittingdirection}.

% Inductive part of the proof of Lemma \ref{splittingdirectionextension}: Induction on J

So we can assume that no vertex $v\in J-\{a',b'\}$ has an edge outside $J$. 
Then $d_J(v)\geq d$ for all $v\in J-\{a',b'\}$ and as $|J|<|G|$ we can apply 
Theorem \ref{conn} to $J$ to find a path of length at least $2^{d}-2$ from $a'$ 
to $b'$. 
Moreover, unless $J$ is isomorphic to $Q_d$ with $a'=a$ and $b'=b$ where $a$ 
and $b$ are at even Hamming distance $J$ contains a path of length at least 
$2^{d}-1$ between $a$ and $b$, so we may assume this is the case. Since $G$ 
is not isomorphic to $Q_d$, the graph $G'=G[V(G)-J\cup \{a,b\}]$ is non-empty 
and all $v$ in $G'-\{a,b\}$ have degree at least $d$ in $G'$. 

If $a$ and $b$ both have more than two limbs in $G'$, $G'$ is $2$-connected. Then 
as $|G'|<|G|$ we can apply Theorem \ref{conn} to $G'$. This gives an $a-b$ path 
in $G'$ of length at least $2^d-1$ unless $G'$ is isomorphic to $Q_d$. Now if $G'$ 
was isomorphic to $Q_d$ then $J$ and $G'$ would both contain the subcube containing 
$a$ and $b$, which has at least four points since $a$ and $b$ are at even Hamming 
distance. But from construction $G'$ and $J$ only share $a$ and $b$, so $G'$ is not 
isomorphic to $Q_d$ and therefore contains an $a-b$ path of length at least $2^d-1$, 
a contradiction.

So one of $a$ and $b$ has exactly one limb. Let this be $a$ say. Then $G'=G_{C'}$ 
for some component $C'$ of $H$ as all limbs of $b$ must have an
out-neighbour in $H$. Again we can apply Theorem \ref{niceextension} to $G'$ to
obtain a $2$-connected subgraph $\widetilde{J}$ and vertices $\widetilde{a}\in
G'_a$ and $\widetilde{b}\in G'_b$. As in Lemma \ref{niceextension}(i) for any $v\in
\widetilde{J}-\{ \widetilde{a},\widetilde{b}\}$, $\widetilde{J}$ contains all
neighbours of $v$ in $G_C=G'$. As such $v$ can have no neighbours in $G$ other
than those in $G'$ we have $d_{\widetilde{J}}(v)=d_{G}(v) \geq d$. Theorem
\ref{conn} holds for $\widetilde{J}$ taking $\widetilde{a}$ and
$\widetilde{b}$ in place of $a$ and $b$. This shows that $\widetilde{J}$ contains 
a $\widetilde{a}-\widetilde{b}$ path of length at least $2^d-2$ which extends to
an $a-b$ path in $G'$ of length at least $2^{d}-2$. As above, since $J$ and 
$\widetilde{J}$ cannot both be isomorphic to $Q_{d}$ if $a$ and $b$ are at even 
Hamming distance, one again must contain an $a-b$ path of length at least $2^d-1$.
This contradicts our assumption and proves the Lemma. \hspace{9.8cm} $\square $
\vspace{.5mm}

% Section 7: End of proof - removing the degree assumption

Lastly, we show that the degree condition can again be removed.

% Analogue of Lemma 6.1 remove above.

\begin{lem} 
\label{removeextension}
Let $G$ be a $2$-connected subgraph of $Q_n$ with $a,b\in G$  
such that $d_G(a)=2$ and $d_G(v)\geq d$ for all $v\in V(G)-\{a,b\}$, 
where $d\geq 3$.
Suppose that Theorem \ref{conn} holds for smaller degrees and 
for all graphs $G'$ with $|G'|<|G|$. 
Furthermore, suppose that $G$ does not contain an $a-b$ path 
of length at least $2^d-1$. 
Then the following hold:
\begin{itemize}
{
\item[\emph{(i)}] $G-a$ is a $2$-connected graph.

\item[\emph{(ii)}] There is a splitting direction $i$ for $G_a$ and $G_b$ 
so that $d_{G_a}(a)\geq 2$ or $d_{G_b}(b)\geq 2$.
}
\end{itemize}
\end{lem}

As the proof is identical to that of Lemma \ref{remove} we will not repeat it. 
We can now finally give the proof of Theorem \ref{conn}.
\vspace{2mm}

% Proof of Main Theorem - 2^d-1 without degree assumption

\noindent \emph{Proof of Theorem \ref{conn}.} The proof is the same as that of Theorem 
\ref{connmod} up until the construction of the interaction digraph $H$. Let $a,a',b, v$ 
and $E_v$ be as in the proof of Theorem \ref{connmod}.

We again take $H$ to be the bipartite multidigraph $H=(A',B,\overrightarrow{E})$ where 
$A'=\{K_1,\ldots ,K_r\}$ and $B=\{L_1,\ldots , L_s\}$, the set of limbs of $a'$ in $G_a$ 
and $b$ in $G_b$ respectively. We also adjoin $\mbox{Core}(b)$ to $B$ if it is non-empty. 
Now from Lemma \ref{listextension}(i), for each endblock $E$ in a limb $K$ of $a$, 
$K\neq \{a,a'\}$, $E$ contains an exit vertex $x_E$ such that $p(x_E)\neq b$ and $E$ 
contains a path of length at least $2^{d-1}-1$ from $\mbox{cutv}(E)$ to $x_E$. Pick one 
such exit vertex $x_E$ for each endblock $E$ other than $\{a,a'\}$ of $G_a$ and such an 
exit vertex $y_F$ for each endblock $F$ of $G_b$. Also let $x_{E_v}=v$. 
For each endblock $E$ of $G_a$, $E\neq \{a,a'\}$, adjoin a directed edge to $H$ from 
$K\in A'$ to $L\in B$ if $E$ is an endblock of $K$ and $p(x_E)\in L$. Similarly, for 
each endblock $F$ of $G_b$, adjoin a directed edge to $H$ from $L\in B$ to $K\in A'$ 
if $F$ is an endblock of $L$ and $p(y_F)\in K$. Again every limb other than $\{a,a'\}$ 
has an outneighbour in $H$.

We now claim that we have a stronger analogue of Propostion \ref{intextension} in this 
case, namely:
\vspace{2mm}

% Section 7: Proof Conn Claim 3

\noindent \textbf{Claim 3:} $p(x_E)\notin \mbox{int}(F)$ for all endblocks $E$ and $F$, $E\neq E_v$. 
\vspace{2mm}

Indeed, if this were the case, then $E$ would contain a $\mbox{cutv}(E)-x_E$ path $P_1$ of 
length at least $2^{d-1}-1$ and $F$ would contain a $p(x_E)-\mbox{cutv}(F)$ path $P_2$ of 
length at least $2^{d-1}-2$. Combining $P_1$ and $P_2$ with the $x_Ep(x_E)$ and extending 
to $a'$ and $b$ we have an $a'-b$ path of length at least $2^d-2$. Appending the edge 
$aa'$ to this, $G$ contains an $a-b$ path of length at least $2^d-1$, a contradiction.

This claim allows us to establish Lemma \ref{path3extension} and \ref{bodyextension} with the 
same proofs as in Lemma \ref{paths3} and \ref{body}. Using this we again find a component $C$ of $H$ 
consisting entirely of limbs and not containing $\{a,a'\}$.

Now apply Lemma \ref{niceextension} again to $G_C$ to find $J$ with properties (i)-(iii). 
As $a'$ is the neigbhour of $a$, we will write $a_J$ and $b_J$ for the vertices of $J$ 
guaranteed by Lemma \ref{niceextension}. In this $J$ we can actually always 
guarantee that both of the paths in (ii) and (iii) have length at least $2^{d-1}$. To see 
this, note that we can replace Claim 1 in the proof of \ref{nice} with the following:
\vspace{2mm}

\noindent \textbf{Claim 4:} $J-b$ has a $2$-connected subgraph $J'$ containing all of $J\cap G_a$ 
and an endblock $F$ of $G_b$.
\vspace{2mm}

\noindent This claim is immediate from the proof of Lemma \ref{nice} on noticing that, by Claim 3, 
$p(S_i)$ cannot lie entirely in the interior of an endblock of $G_b$ and so $\mbox{span} _{G_b}(S_i)-\{b\}$ must 
contain an endblock of $G_b$.

Now use $J'$ as in the final paragraph of the proof of Lemma \ref{nice} and choose one of $E$ 
or $F$ in place of $E$ so that $v$ is on the opposite side of $J'$ to the chosen endblock; 
choose $E$ if $v\in G_b$ and $F$ if $v\in G_a$. The path then constructed contains an endblock path 
of length at least $2^{d-1}-2$ and a path of length at least $2$ from $v$ to the chosen 
endblock. This gives the claimed path of length at least $2^{d-1}$.

We now use this stronger fact to complete the proof. If there are no edges from $v\in J-\{a_J,b_J\}$ 
to a vertex in $G-J$ then, by induction on Theorem \ref{conn}, $J$ contains an $a_J-b_J$ path of length 
at least $2^d-2$. Extending this to an $a'-b$ path and appending the edge $aa'$ we obtain an $a-b$ path 
of length at least $2^d-1$, a contradiction.

So we can assume that some $v\in J-\{a_J,b_J\}$ has a neighbour $w$ outside $J$. The proof can now be 
finished in exactly the same way as the proof of Lemma \ref{splittingdirection}. As here we always 
adjoin one of the paths in $J$ of length at least $2^{d-1}$ with another endblock path of length at 
least $2^{d-1}-2$ via the edge $vw$, the $a-b$ path we create always has length at least $2^d-1$. 
This contradiction proves the Theorem. \hspace{10cm} $\square$

\section{Generalizations}

The reader might have noticed that we have used very little about $Q_n$ in the proof
of Theorem \ref{conn}. The $n$-dimensional grid ${\Z}^n$ is the graph whose
vertex set consists of $n$-tuples with entries in $\Z$ and in which two vertices
$x$ and $y$ are adjacent if $|x_i-y_i|=1$ for some $i\in [n]$ and $x_j=y_j$ for
all $j\neq i$. The next theorem extends Theorem \ref{conn} (and therefore
Theorems \ref{mainthm} and \ref{cyclethm}) to subgraphs of $\Z^n$.

\begin{thm}
{
Let $G$ be a $2$-connected subgraph of ${\Z}^n$ and $a,b\in V(G)$. Suppose that
$d(z) \geq d$ for all $z \in V(G)-\lbrace a,b\rbrace$. Then $a$ and $b$ are
joined by a path of length at least $2^d-2$. Furthermore unless $G$ is
isomorphic to $Q_d$ with $a$ and $b$ at even Hamming distance from each other,
$G$ contains an $a-b$ path of length $2^d-1$.
}
\end{thm}

\begin{proof}
{
The crucial property of $\Z ^n$ here is that we can always find a splitting of
$G$ into two connected pieces, $G_a$ and $G_b$ with $a\in G_a$ and $b\in G_b$
such that $d_{G_a}(a)\geq 1$ and $d_{G_b}(b)\geq 1$ and all $v\in G$ lose at
most one neighbour in their piece. Indeed, taking some coordinate $j$ on which
$a$ and $b$ differ, say with $a_j>b_j$, let $G_1$ be the induced subgraph of $G$
consisting all vertices $v$ with $v_j\geq a_j$ and $G_2$ be the induced subgraph
of $G$ consisting of all $w$ for which $w_j<a_j$. Again with the same
modification to these graphs as in Lemma \ref{splitting} we obtain connected
graphs $G_a$ and $G_b$ with the required degree conditions. From here on the
proof is identical to that of Theorem \ref{conn}.
}
\end{proof}

Moreover, the same proof also extends to subgraphs of the discrete torus
$C_k^{n}$ provided $k\geq 4$. Now we cannot expect a bound of the form $C2^d$ as
above for subgraphs of the discrete torus $C_3^{d}$ as this graph has minimum
degree $2d$ but only $3^d$ points. This shows that given a subgraph $G$ of
$C_3^{n}$ of minimal degree at least $d$ we cannot in general guarantee a path
of length more than $3^{\frac{d}{2}}-1$ in $G$. 

Why does our approach not work in this case? The main reason is that we cannot
guarantee a partition into two subgraphs such that all vertices lose at most one
neighbour in their piece. Can we still guarantee an exponentially long path in
this case? 

The following general result shows that we can. 

\begin{thm}
{
\label{gen} Let $k\in \N$ and $G$ be a $2$-connected graph with $a,b\in V(G)$.
Suppose $d(v)\geq d$ for all $v\in V(G)-\{a,b\}$. Furthermore, suppose that $G$
has the following property:
\begin{quotation}
{\noindent Given any two vertices $x,y\in G$, there is a partition of $V(G)$
into two sets $X$ and $Y$ 
with $x\in X$ and $y\in Y$ such that $d_{G[X]}(v)\geq d(v)-k$ for all $v\in X$
and $d_{G[Y]}(v)\geq d(v)-k$ for all $v\in Y$.}
\end{quotation}
Then $G$ contains an $a-b$ path of length at least $2^{\frac{d}{k+2}}$.
}
\end{thm}

Note that if the property above holds for $G$, it also holds for all subgraphs
of $G$. Also note that Theorem \ref{extension} immediately follows from Theorem 
\ref{gen}. As an immediate corollary of Theorem \ref{gen} we have the following:

\begin{cor}
 {
 Every subgraph of $C_3^n$ of minimum degree at least $d$ contains a path of
length at least $2^{\frac{d}{4}}$.
 }
\end{cor}

It would be interesting to decided what the correct lower bounds for the length
of the longest path in subgraphs of $C_3^n$ with minimum degree at least $d$.

\begin{conjecture}
 {
 Given a subgraph $G$ of $C_3^n$ with minimum degree at least $d$, $G$ must
contain a path of length at least $3^{\frac{d}{2}}-1$.
 }
\end{conjecture}

Another consequence of Theorem \ref{gen} is the following result for product
graphs.

\begin{thm}
{
Let $G_1,\ldots ,G_l$ be graphs with maximum degree at most $k$. Then given any
subgraph $G$ of the Cartesian product graph $\prod_{i=1}^{l}G_i$  of minimum
degree at least $d$, $G$ contains a path of length at least $2^{\frac{d}{k+2}}$.
}
\end{thm}

The proof of Theorem \ref{gen} is similar to that of Theorem \ref{conn} but
shorter.

\begin{proof}
{
The proof is again by induction on $d$. It suffices to prove the result for
$d\geq k+4$ as otherwise it follows from 2-connectivity. As in the proof of
Theorem \ref{conn} we wish to split $G$ into two subgraphs $G_a$ and $G_b$ with
$a\in G_a$ and $b\in G_b$, which is the motivation for the above splitting
property. However, simply taking $a$ and $b$ in place of $x$ and $y$ might not
be useful as both $a$ and $b$ can have degree as low as two in $G$ in which case
in the partition guaranteed above $a$ may end up with all its neighbours in
$Y$. 
Instead we pick a neighbour $a'\neq b$ of $a$ and a neighbour $b'\neq a$ of $b$.
The fact that $G$ is $2$-connected ensures it is possible to pick $a'\neq b'$.
Now take the partition guaranteed from our splitting property above with $x=a'$
and $y=b'$. Moving $a$ to $X$ and $b$ to $Y$ as needed we have that
$d_{G[X]}(v)\geq d(v)-k-1$ for all $v\in X-a$ and similarly for $v\in Y-b$. Both
$a$ and $b$ now have at least 1 neighbour in $G[X]$ and $G[Y]$ respectively.
Finally, denoting the connected component of $G[X]$ containing $a$ by $C_a$, let
$G_b$ be the connected component of $G-C_a$ containing $b$ and
$G_a=G[V(G)-V(G_b)]$. Note that $G_a$ and $G_b$ are connected with $a\in G_a$,
$b\in G_b$. Moreover, $d_{G_a}(v)\geq d(v)-k-1$ for all $v\in G_a-a$ and
$d_{G_b}(v)\geq d(v)-k-1$  for $v\in G_b-b$.

We will again analyse the block-cutvertex decompositions of $G_a$ and $G_b$. The
following lemma will be very useful below.

\begin{lem}
{
\label{endblock}
Let $E$ be an endblock of $G_a$ or $G_b$ with $a,b\notin \mbox{int}(E)$. Then
given any two vertices $u,v\in \mbox{int}(E)$, $G[E]$ contains a path of length
at least $2^{\frac{d-k-2}{k+2}}$ from $u$ to $v$.
}
\end{lem}

\begin{figure}
\centering
\input{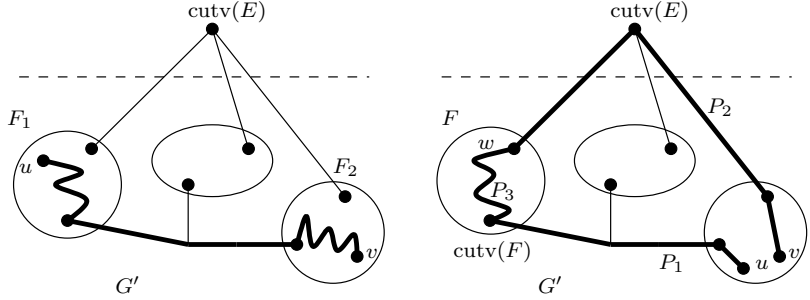}
\caption{Cases where $G'$ is not $2$-connected in Lemma \ref{endblock}}
\label{smallerlongpath}
\end{figure}

\begin{proof}
{
Look at the block-cutvertex decomposition of $G'=G[E]-\mbox{cutv}(E)$. Since $E$
is $2$-connected (as $d\geq k+4$), $G'$ is connected and $\mbox{cutv}(E)$ must
have a neighbour in the interior of every endblock of $G'$. Note that every
vertex $v\in G'$ has $d_{G'}(v) \geq d_G(v)-k-2$. In particular, since $d\geq
k+4$ each endblock $F$ of $G'$ is $2$-connected and has at least three vertices
so that we can by induction apply Theorem \ref{gen} to it. If $G'$ is
$2$-connected then by induction on Theorem \ref{gen} $G'$ contains the desired
path from $u$ to $v$. Thus we may assume that $G'$ is not $2$-connected. If
$u\in \mbox{int}(F_1)$ and $v\in \mbox{int}(F_2)$ where $F_1$ and $F_2$ are two
distinct endblocks of $G'$ then by induction on Theorem \ref{gen}, $G[F_1]$ and
$G[F_2]$ contain $u-\mbox{cutv}(F_1)$ and $\mbox{cutv}(F_2)-v$ paths
respectively, each of length at least $2^{\frac{d-k-2}{k+2}}$. Joining
$\mbox{cutv}(F_1)$ to $\mbox{cutv}(F_2)$ by a third path in $G'$ and combining
all three of these paths, we get a $u-v$ path of length at least
$2^{\frac{d}{k+2}}$, as required. Therefore since $G'$ contains at least two
endblocks, we can assume that one of these, say $F$, does not contain $u$ or $v$
in its interior. Contracting $\mbox{int}(F)$ down to a single vertex in $G[E]$,
the resulting graph is still $2$-connected. Therefore, as in the proof of Lemma
\ref{nice}, $G[E]$ contains two vertex disjoint paths $P_1$ and $P_2$ from the
set $\{u,v\}$ to $\{\mbox{cutv}(F),w\}$ for some $w\in \mbox{int}(F)$, with
$(P_1\cup P_2)\cap (F-\{\mbox{cutv}(F),w\})=\emptyset$. Now using induction on
Theorem \ref{gen} in $F$, it contains a path $P_3$ of length
$2^{\frac{d-k-2}{k+2}}$ from $\mbox{cutv}(F)$ to $w$. Piecing $P_1$, $P_2$ and
$P_3$ together we obtain our desired path.
}
\end{proof}

Again we have:

\begin{prop}
{
\label{intnew}
Let $E$ be an endblock of $G_a$ not containing $a$ and $F$ an endblock of $G_b$
not containing $b$. Then $G$ does not contain an edge from $\mbox{int}(E)$ to
$\mbox{int}(F)$
}
\end{prop}

\begin{proof}
{
Exactly as in Proposition \ref{int}.
}
\end{proof}

\begin{lem}
{
\label{newlist}
We have the following:
\begin{itemize}
{
\item[\emph{(i)}] Given any endblock $E$ of $G_a$ not containing $a$, there are
two disjoint edges from $\mbox{int}(E)$ to $G_b$ in $G$.

\item[\emph{(ii)}] $G_a$ contains an endblock not containing $a$.
}
\end{itemize}
}
\end{lem}

\begin{proof}
{
(i) $E$ must have an exit vertex $x_1$, with neighbour $y\in G_b$, as $G$ is
$2$-connected. If it had only one, $G'=G[E]$ is $2$-connected and every $v\in
G'-\{x_1,\mbox{cutv}(E)\}$ has degree at least $d$. Therefore, by induction on
Theorem \ref{gen}, $G'$ contains a path of length at least $2^{\frac{d}{k+2}}$
from $x_1$ to $\mbox{cutv}(E)$. Extending this path from $\mbox{cutv}(E)$ to $a$
in $G_a$ and from $y$ to $b$ in $G_b$ we obtain an $a-b$ path of desired length.
Therefore we may assume $E$ contains a second exit vertex $x_2$. Now if the
vertices in ${\mbox{int}(E)}$ were only adjacent to $y$ in $G_b$, $x_1y$ and
$x_2y$ must be edges of $G$. Then $G''=G[E\cup \{y\}]$ is $2$-connected and
$d_{G''}(v)\geq d$ for every $v\in G''-\{y,\mbox{cutv}(E)\}$. By Theorem
\ref{gen} $G''$ contains a $\mbox{cutv}(E)-y$ path of length at least
$2^{\frac{d}{k+2}}$. Again, extending this to a path from $a$ to $b$, we have an
$a-b$ path of length at least $2^{\frac{d}{k+2}}$. Therefore we may assume the
two edges exist or we are done.

(ii) The proof is almost identical to the proof of Lemma \ref{list}(ii).
}
\end{proof}

Take an endblock $E$ of $G_a$ not containing $a$, as guaranteed by Lemma
\ref{newlist}(ii). We can choose $E$ such that $a$ and all $v\in G_a-E$ not
contained in the interior of an endblock of $G_a$ lie in the same connected
component of $G_a-E$ (e.g. pick a block $B$ in $G_a$ containing $a$ and choose
$E$ to be a block at maximum distance from $B$ in $\mathcal{B}(G_a)$). Let
$x_1y_1$ and $x_2y_2$ be the disjoint edges of $G$ with $x_1,x_2\in
\mbox{int}(E)$ and $y_1,y_2\in G_b$ guaranteed by Lemma \ref{newlist}(i). By
Proposition \ref{intnew}, $y_1,y_2\notin \mbox{int}(F)$ for all endblocks $F$ of
$G_b$ not containing $b$ in its interior. 

\begin{figure}
\centering
\input{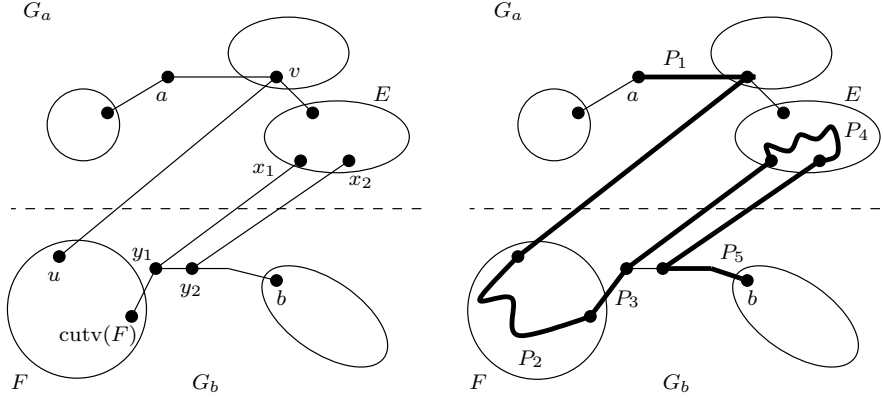}
\caption{Path created in Theorem \ref{gen}}
\label{Generalized}
\end{figure}

Now looking at the block-cutvertex decomposition of $G_b$ we can choose two
vertex disjoint paths in $G_b$ from $\{y_1,y_2\}$ to $\{b, \mbox{cutv}(F)\}$
where $F$ is some endblock of $G_b$ not containing $b$. Lets say that these
paths are $P_3$ from $\mbox{cutv}(F)$ to $y_1$ and $P_5$ from $y_2$ to $b$.
Applying Lemma \ref{newlist}(i) to $F$ we see that there exists $u\in
\mbox{int}(F)$ adjacent to some $v\in G_a$, $v\neq \mbox{cutv}(E)$. Furthermore,
by Proposition \ref{intnew} $v\notin \mbox{int}(E')$ for any endblock $E'$ of
$G_a$. From our choice of $E$ there exists an $a-v$ path $P_1$ in $G_a-E$.
Finally by induction on Theorem \ref{gen}, $F$ contains a $u-\mbox{cutv}(F)$
path $P_2$ of length at least $2^{\frac{d-k-1}{k+2}}$ and by Lemma
\ref{endblock} $E$ contains an $x_1x_2$ path $P_4$ of length at least
$2^{\frac{d-k-2}{k+2}}$. Combining these five paths we obtain an $a-b$ path
$P=P_1vuP_2P_3y_1x_1P_4x_2y_2P_5$ of length at least
$2^{\frac{d-k-1}{k+2}}+2^{\frac{d-k-2}{k+2}}>2^{\frac{d}{k+2}}$ as required. 
}
\end{proof}

The cycle analogues of the above theorems can be obtained in a similar fashion
to the proof of Theorem \ref{cyclethm} from Theorem \ref{conn}.

As mentioned in the Introduction, we do not know the correct bound for the
length of the longest path in a subgraph of $Q_n$ when the minimum degree
condition in Theorem \ref{mainthm} is replaced by an average degree condition.
Is the following possible?
\vspace{0.25cm}

\begin{conjecture} 
 {
Every subgraph of $Q_n$ with average degree at least $d$ contains a path of
length at least $2^d-1$.
 }
\end{conjecture}

\section*{Acknowledgements}
{
I would like to thank my supervisor Imre Leader for many helpful conversations. I would also 
like to thank David Conlon for reading through an earlier version of this paper and the 
anonymous referee for his/her valuable comments.
}


\begin{thebibliography}{99}
\bibitem{bernstein} A.J. Bernstein: Maximally connected arrays on the $n$-cube,
\textit{SIAM J. Appl. Math.} \textbf{15}(1967), 1485-1489.
\bibitem{com} B. Bollob\'as: \textbf{Combinatorics: Set Systems, Hypergraphs,
Families of Vectors and Combinatorial Probability}, Cambridge University Press,
1st ed, 1986.
\bibitem{mgt} B. Bollob\'as: \textbf{Modern Graph Theory}, Springer, 1st ed,
1998.
\bibitem{dirac} G.A. Dirac: Some theorems on abstract graphs, \textit{Proc.
London Math. Soc.} \textbf{2} (1952), 69-81.
\bibitem{harper} L.H. Harper: Optimal assignments of numbers to vertices,
\textit{SIAM J. Appl. Math.} \textbf{12}(1964), 131-135.
\bibitem{hart} S. Hart: A note on edges of the $n$-cube, \textit{Discrete Math.}
\textbf{14}(1976), 157-163.
\bibitem{lindsey} J.H. Lindsey: Assignment of numbers to vertices, \textit{Amer.
Math. Monthly} \textbf{71}(1964), 508-516.
\end{thebibliography}
\end{document}